\documentclass[final]{siamart1116}
\usepackage{graphicx,helvet}
\usepackage{epstopdf}
\usepackage[percent]{overpic}
\usepackage{xcolor,amssymb}

\usepackage{epsfig}
\usepackage{cite}
\usepackage[english]{babel}
\usepackage[utf8]{inputenc}

\allowdisplaybreaks

\newcommand{\eps}{\varepsilon}
\renewcommand{\epsilon}{\varepsilon}

\newcommand{\R}{\mathbb{R}}

\newcommand{\norm}[1]{\left\|#1\right\|}

\newcounter{dctr}[section]
\numberwithin{equation}{section}

\newtheorem{remark}{Remark}[section]

\newcommand{\bigO}{\mathcal{O}}

\renewcommand{\L}{\mathcal{L}}
\def\sii{\Leftrightarrow}
\def\Jr{\{0,\dots,r-1\}}

\def\mz{\alpha}
\def\mo{\alpha+1}

\newcommand{\N}{\mathbb{N}}

\def\sf#1#2{\sum_{i} \frac{#1_i #2_i}{v_i-\lambda}}

\numberwithin{equation}{section}
\numberwithin{theorem}{section}
\numberwithin{table}{section}
\numberwithin{figure}{section}

\headers{High-order WENO reconstructions}{Baeza, B\"{u}rger, Mulet, and Zor\'{\i}o}

\allowdisplaybreaks
\numberwithin{equation}{section}

\title{WENO reconstructions of \\ unconditionally optimal high order} 

\date{\today} 

\author{Antonio Baeza\thanks{Departament de Matem\`{a}tiques, 
Universitat de Val\`{e}ncia, Av.\ Vicent Andr\'es Estell\'es, E-46100
Burjassot,  
 Spain.  E-Mail: 
   {\tt antonio.baeza@uv.es}}
\and
Raimund  B\"urger\thanks{CI$^{\mathrm{2}}$MA and De\-par\-ta\-mento de
     Ingenier\'ia Matem\'a\-tica, Uni\-ver\-sidad de Con\-cepci\'on, Casilla
     160-C, Con\-cep\-ci\'on, Chile. E-mail: \texttt{rburger@ing-mat.udec.cl}}
     \and  
     Pep Mulet\thanks{Departament de Matem\`{a}tiques, 
Universitat de Val\`{e}ncia, Av.\ Vicent Andr\'es Estell\'es, E-46100
Burjassot,  
 Spain.  E-Mail: 
   {\tt pep.mulet@uv.es}}
\and 
  David Zor\'{\i}o\thanks{CI$^{\mathrm{2}}$MA,   Universidad de
Concepci\'{o}n, Casilla 160-C, Concepci\'{o}n, Chile.  E-Mail:
  {\tt dzorio@ci2ma.udec.cl}}
  }

\begin{document}
\maketitle 

\begin{abstract}
A modified Weighted Essentially Non-Oscillatory (WENO) reconstruction technique preventing accuracy loss near critical points 
 (regardless of their order) of the underlying data  is presented. This approach only  uses  local data from the reconstruction stencil 
  and does not rely on any sort of scaling parameters. The key novel
  ingredient is a weight design based on a new smoothness indicator, which defines the first WENO reconstruction procedure  that never loses accuracy on smooth data, regardless of the presence of critical points  of any order, and is therefore addressed as optimal WENO (OWENO) method.  
   The corresponding weights  are non-dimensional and scale-independent. 
 The weight designs are supported by  theoretical results  concerning  the accuracy  of the smoothness indicators.
  The method is validated by numerical tests related to 
    algebraic equations, 
scalar conservation laws,  and systems of conservation laws.
\end{abstract} 



\section{Introduction}

\subsection{Scope}

\label{subsec:scope} Weighted Essentially Non-Oscillatory (WENO) reconstructions, initially proposed by 
 Liu  et al.\ \cite{LiuOsherChan94} and later improved by  Jiang and  Shu \cite{JiangShu96}, have become a common  ingredient  of high-resolution schemes for the 
 numerical solution of  hyperbolic conservation laws. The standard  initial-value 
 problem  is of the type 
 \begin{align}  \label{goveq} 
  \boldsymbol{u}_t + \sum_{i=1}^{\mathcal{D}} \boldsymbol{f}_i ( \boldsymbol{u} )_{x_i} 
   = \boldsymbol{0},  \quad \boldsymbol{x} = (x_1, \dots,x_{\mathcal{D}}) \in \mathbb{R}^{\mathcal{D}},  \quad t>0, 
  \end{align}  
where 
    $\boldsymbol{u} = \boldsymbol{u} (\boldsymbol{x}, t) = (u_1, \dots, u_N)^{\mathrm{T}}$ 
     is the  vector of sought unknowns and $\boldsymbol{f}_i ( 
     \boldsymbol{u} ) = (f_{i,1}(\boldsymbol{u} ), \dots, f_{i,N} (  \boldsymbol{u} ))^{\mathrm{T}}$  are given flux vectors, supplied with an initial condition
    \begin{align} \label{initcond} 
     \boldsymbol{u} (\boldsymbol{x}, 0) = \boldsymbol{u}_0 ( \boldsymbol{x}), \quad 
      \boldsymbol{x} \in \mathbb{R}^{\mathcal{D}}. 
      \end{align} 
   Such schemes  (in short, ``WENO schemes'')  present a high order of accuracy in smooth zones and,
   through a sophisticated construction of non-linear weights~\cite{JiangShu96}, 
    avoid the oscillatory behaviour 
typical of the reconstructions from discontinuous data. 
 However, such weights are sensitive not only to discontinuities, but also  to abrupt changes in any higher derivative of the function that generates the data, which leads to an undesired loss of accuracy near critical points. A variety of solutions to handle this problem have been proposed; see for instance    \cite{SINUM2011,ArandigaMartiMulet2014,HenrickAslamPowers2005,YamaleevCarpenter2009}. 
  However, none of them allows to unconditionally attain  the optimal order of accuracy (that is, regardless of the order of the critical points) depending only on the local data without ending up with dimensional or scale-dependent weights. In other words,  either  some dimensional (namely,  
   grid-size-dependent) or not properly scaled parameter is used, or data from the global numerical solution are employed to define a non-dimensional and scale-independent parameter to prevent such loss of accuracy.

It is the purpose of  this paper
  to  design weights in such a way  that the    associated reconstruction algorithm does not lose accuracy in smooth zones, even in presence of critical points of any order. The decisive novelty of the new  non-dimensional and scale-independent weights that only  use  information from the local data of  the stencil. Since the order of accuracy of the resulting new WENO schemes is optimal, we refer to them 
    as ``optimal WENO'' (OWENO) schemes. 
    At the core of this paper is an analysis of the 
   accuracy properties involving the asymptotics  of the smoothness indicators as the grid size goes to zero. 
    This issue   has often  been studied  for reconstructions of specific orders along the literature, but no full proof for the general case 
     has been advanced so far. We provide such a proof. The theoretical tools will be then available to fully and solidly analyze the accuracy of the reconstructions proposed, and are utilized to design OWENO reconstructions of 
       unconditionally optimal order of accuracy regardless of the order of the critical points.
\subsection{Related work}
Overviews on WENO schemes include  \cite{shu98,shu09,zhang16}. The particular  problem of  achieving optimal order of accuracy near critical points is tackled in many  works. Henrick  et al.\ \cite{HenrickAslamPowers2005}    obtain optimal order convergence near critical points for the case of fifth order through a simple modification of the  weights by  Jiang and Shu \cite{JiangShu96}, which involves mapping the weights to values that  satisfy an optimality condition. The approach was further extended up to order~17 by  Gerolymos 
  et al.\  \cite{GerolymosSenechalVallet2009}   and further enhanced by    Feng et al.\ \cite{FengHuWang2012}   by means of a different mapping.
 A different weight design was followed in the fifth-order WENO-Z method by  Borges  et al.\ \cite{BorgesCarmonaCostaEtAl2008}, which attains fourth-order accuracy even at critical points. Castro et al.\ \cite{CastroCostaDon2011}    
   extended the WENO-Z scheme  to any odd order of accuracy, achieving optimal order at critical points by proper parameter tuning. 

Following an idea  similar to that of  WENO-Z schemes,  Yamaleev and Carpenter \cite{YamaleevCarpenter2009}  introduced    a new method, named ESWENO, based on the third-order case previously introduced  in~\cite{YamaleevCarpenter2009a}  that ensures energy stability in an $L^2$ norm. Even though it was not their primary goal to enhance order at critical points, it turns out that the resulting scheme achieves optimal order in the presence of critical points provided that the number of zero derivatives is at most the order of the scheme minus three. 
Another way to handle the problem of the order loss at critical points is the modification of the smoothness indicators.  Ha et al.\ \cite{HaKimLee2013} 
proposed a new smoothness measurement that  provides optimal order for functions with
critical points, but in which the second derivative is not zero.

The design of weights in WENO schemes typically involves a quantity $\varepsilon$ that avoids division by zero whenever an smoothness indicator becomes zero. This parameter  was set to a fixed quantity $\varepsilon=10^{-6}$ in \cite{JiangShu96}, 
	but Ar\`{a}ndiga   et al.\ \cite{SINUM2011}  noted that the
        choice of $\varepsilon$ is 
crucial for the achievement of optimal order at critical points and that, for the case of 
the original weights of Jiang and Shu \cite{JiangShu96}, 
 the choice of $\varepsilon$ proportional to the square of the mesh size provides
 the desired accuracy even at critical points.  A similar analysis was later performed 
 by  Don and  Borges \cite{DonBorges2013}, regarding WENO-Z
schemes, and  in \cite{ArandigaMartiMulet2014}  with respect to the ESWENO weights of Yamaleev and Carpenter (see also \cite{Kolb2014}), thus requiring a 
scale-dependent parameter.

\subsection{Outline of this paper}

The required theoretical
background of this work is outlined  in Section~\ref{prel}, where we derive bounds involving the asymptotical behaviour of the smoothness indicators as the grid size tends to zero.  After collecting some preliminaries of notation in Subsection~\ref{subsec:prel}, 
      we state in Subsection~\ref{subsec:wenorec} some results that will be helpful for the analysis of the accuracy of  WENO reconstructions for both cases of point value and cell average data.  
Such bounds are   the  key ingredients within
Section~\ref{wd}, which is devoted to the definition of the
new OWENO  reconstructions that attain optimal order of accuracy
regardless of the number of consecutive zero derivatives of the
function to be reconstructed, and without using any scaling
parameter. Inside this section, we first motivate the issues involving
the accuracy loss of the existing schemes in the literature in
Subsection \ref{subsec:mot}. Then, we propose a novel smoothness
indicator which overcomes these issues in Subsection~\ref{subsec:nsi},
which is the main novelty of this paper, along with 
theoretical results that support the considerations on the
optimal accuracy. Finally, Subsection~\ref{sec:alg1} 
summarizes the algorithm of the proposed method with the novel
smoothness indicators.

In summary, we prove that the new scheme has unconditionally optimal order of accuracy under 
those conditions, therefore overcoming the issue of the scheme proposed by Yamaleev  and Carpenter \cite{YamaleevCarpenter2009}  involving the accuracy loss near critical points in which the number of consecutive vanishing derivatives is the order of the scheme minus two.  In Section~\ref{sec:num_exp} we present some
numerical experiments, both for algebraic problems in Subsection
\ref{subsec:alg_prob} and problems involving hyperbolic conservation
laws in Subsection \ref{subsec:cons_laws}. Finally, in Section
\ref{sec:conclusions} some conclusions are drawn. Some  technical results related 
 the accuracy of  OWENO schemes are collected in Appendix~\ref{sec:appendix}.

\section{Regularity properties of functions and smoothness indicators}\label{prel}

The analysis of WENO schemes will be carried out in one space dimension, where~$x$ denotes the spatial coordinate and $h>0$ is the uniform 
 mesh width.  
This section is devoted to analyze the asymptotic accuracy
  properties of the smoothness indicators  by Jiang and Shu
  \cite{JiangShu96}, which on a stencil of $2r-1$ points
  $\{x_{-r+1},\ldots,x_{r-1}\}$, with $x_{i+1}=x_i+h$, have the form
\begin{equation*}
  I=\sum_{l=1}^{r-1}\int_{x_0-h/2}^{x_0+h/2}h^{2l-1} \bigl(p^{(l)}(x)\bigr)^2 
   \, \mathrm{d}x,
\end{equation*}
where $p$~is a reconstruction polynomial,
    corresponding to a substencil of $r$ points.
The key  result  that lays the foundation  for the ulterior accuracy analysis of the new OWENO reconstructions  is Theorem~\ref{is}, stated below, which provides the {\it exact} convergence rate of the Jiang-Shu smoothness indicators near critical points of any order. In order to prove this result, some technical definitions and results will be presented before.
 Theorem~\ref{is}  is also crucial for the accuracy analysis near critical points of all the WENO reconstructions modalities presented in the literature 
  that are based on  the Jiang-Shu smoothness indicators.

\subsection{Preliminaries}  \label{subsec:prel} 
 For a piecewise smooth function  with jump
discontinuities $f:\R \to\R$,  we
use the standard  notation $f(h)=\bigO(h^\alpha)$ for $\alpha 
 \in \mathbb{Z}$ to indicate 
 the behaviour of a function~$f$ as $h \to 0$ in the standard sense, that is,  
\begin{align*}
  f(h)&=\bigO(h^\alpha)\sii \limsup_{h\to 0} \bigl|f(h) h^{-\alpha}\bigr| <\infty.
 \end{align*} 
 Furthermore, we write $f(h) =\bar{\bigO}(h^\alpha)$ to express the more restrictive property
 \begin{align*} 
  f(h)&=\bar{\bigO}(h^\alpha)\sii \limsup_{h\to 0} \bigl|f(h) h^{-\alpha}\bigr|  <\infty \quad \text{and} \quad \liminf_{h\to 0} \bigl|
    f(h) h^{-\alpha}\bigr| > 0.
\end{align*}
It follows  for $\alpha, \beta \in \mathbb{Z}$   that
$\bar\bigO(h^{\alpha})^{-1} = \bar\bigO(h^{-\alpha})$,
$\bigO(h^{\alpha}) \bigO(h^{\beta}) = \bigO(h^{\alpha+\beta})$ 
and   $\smash{\bar\bigO(h^{\alpha}) \bar\bigO(h^{\beta}) =
\bar\bigO(h^{\alpha+\beta})}$. Moreover, we say that a function $f$ has a  critical point of 
order~$k\geq0$ at $x$ if $\smash{f^{(l)}}(x)=0$ for 
$l=1,\dots,k$ and $\smash{f^{(k+1)}}(x)\neq 0$. For $k=0$ this includes 
 the degenerate case of a point~$x$ at which $f'(x) \neq 0$.

We extend the classical notation for continuously higher
differentiable function to denote by $f\in  C^{s}(z)$ if
there exists $\delta > 0$ such that 
$f\in  C^{s}(z-\delta, z+\delta)$ and by $f\in
 C^{s}(z^{\pm})$ if there exists $\delta > 0$ such that $f$~is $s$~times continuously differentiable in $(z-\delta,
z+\delta)\setminus \{z\}$ and  $\smash{\lim_{x\to z^{\pm}}
  f^{(s)}(x)=f^{(s)}(z)}$.

\subsection{WENO reconstructions} \label{subsec:wenorec}

 For a  stencil
 \begin{align} \label{Sdef} 
 S=\{x_{-r+1},\ldots,x_{r-1}\}
 \end{align} 
 of $2r-1$ points~$x_j=x_{j,h}$, where $x_{j+1}-x_j=h$ for $-r+1\leq j\leq r-1$,  
 and a  scalar function~$f$ we assume that the data
 $\{f_{-r+1},\ldots,f_{r-1}\}$ are  either 
 point values \begin{align} \label{fpointvalues} 
f_j=f(x_j), \quad -r+1\leq j\leq r-1,
\end{align} 
 or   cell averages 
 \begin{align} \label{fcellaverages} 
 f_j=\frac{1}{h}\int_{x_{j-1/2}}^{x_{j+1/2}}f(x) \, \mathrm{d}x, \quad -r+1\leq j\leq r-1,
 \end{align} 
where in both cases we wish to approximate the point value~$f(x_{1/2})$. 

 We  denote by $\Pi_k$, $k \in \mathbb{N}_0$, the space of polynomials of {\em maximal} degree~$k$, and by~$\smash{\bar{\Pi}_k}$ the space of polynomials of {\em exact} degree~$k \in \mathbb{N}_0$.  
 Let  $p_{r,i} \in \Pi_{r-1} $
 denote the reconstruction polynomial  of the
 substencils 
 \begin{align} \label{Sridef} 
 S_{r,i}=\{x_{-r+1+i},\ldots,x_i\},  \quad 0\leq
i\leq r-1, 
\end{align} 
with the  
 interpolation property  
 $p_{r,i} ( x_j) = f_j$ for reconstructions from point
   values \eqref{fpointvalues}  or 
   \begin{align*} 
   \int_{x_{j-1/2}}^{x_{j+1/2}}p_{r,i}(x)\, \mathrm{d}x=f_j
   \end{align*} 
    for
   reconstructions from cell averages \eqref{fcellaverages} for all $x_j \in S_{r,i}$. In what follows, we omit  
 the  subindex~$r$ when no confusion may arise.

The WENO  strategy consists in
defining a reconstruction $q$ as a convex combination 
$ q(x_{1/2})= \omega_0 p_{0} (x_{1/2}) + \omega_1 p_1 (x_{1/2}) + \dots + \omega_{r-1} p_{r-1} (x_{1/2})$
of the
individual reconstructions~$p_i$ with appropriately designed weights
$\omega_0, \dots, \omega_{r-1}\geq 0$, where $\omega_0 + \cdots +
\omega_{r-1}=1$, which satisfy $\omega_i\approx c_i$
  on smooth zones, with $c_i$ the linear ideal weights
  \cite[Proposition 2]{SINUM2011}, satisfying that $c_0 p_{0} (x_{1/2}) + c_1
  p_1 (x_{1/2}) + \dots + c_{r-1} p_{r-1} (x_{1/2})$ coincides with the interpolatory
polynomial of order $2r-1$ at $x_{1/2}$. 
The weights $\omega_i$ are functions of some 
smoothness indicators, which we take according to Jiang and Shu \cite{JiangShu96}:
\begin{equation}\label{eq:SI}
  I_{i}= \sum_{l=1}^{r-1}  \int_{x_{-1/2}}^{x_{1/2}}h^{2l-1} \bigl(p_{i}^{(l)}(x)\bigr)^2 
   \, \mathrm{d}x.
\end{equation}
Notice that $I_{i}=0$ implies that 
$\smash{p'_{i}}=0$   on an 
interval of positive length, so that $\smash{p'_{i}}$ is  zero everywhere, i.e., 
$f$~is constant at the points of~$S_{r,i}$. 
  
\begin{theorem}\label{is}
  Let $z,\alpha\in\R$, $h>0$ and $x_i=z+(\alpha+i)h$, $-r+1\leq i\leq r-1$ define  a stencil of equally-spaced nodes. If $f$ has a critical point of order~$k$ at $z$, then 
   the Jiang-Shu smoothness indicator \eqref{eq:SI} satisfies 
    $\smash{I_{i}=\bar\bigO(h^{2\kappa})}$, where  
  \begin{align*}
    \kappa=\begin{cases}
      \min\{l\in\N\colon  2 | l, l\geq k,
      f^{(l+1)}(z)\neq0\} & \text{\em for $r=2$ and   $\alpha+i=1/2$,} \\
      k+1 & \text{\em otherwise.}
    \end{cases}      
  \end{align*}
\end{theorem}

\begin{proof}
  From Lemma \ref{zr} applied to $n=r-1$ and  \eqref{cijdef}, we obtain
  \begin{align}\label{eq:225}
  p_{i}^{(j)}(z+wh)=\sum_{s=j}^mb_{i,s,j}(w)h^{s-j}f^{(s)}(z)+\bigO(h^{m+1-j}),
\end{align}
where $\smash{b_{i,s,j}}$ denotes the function $b_{s,j}$ given by Lemma~\ref{zr} corresponding 
 to the stencil~$S_{r,i}$. Notice that the condition $\alpha+i=1/2$ is
equivalent to $\smash{a_{0,i}=-a_{1,i}}$.
We apply \eqref{eq:225} for $m=\kappa:= \min \{\nu \in \mathbb{N}  \colon b_{i,\nu,1}(w) \smash{f^{(\nu)}(z)} \neq 0
\}$.  
Then by the  definition of $k$ and Lemma~\ref{zr}
   we get for $j\leq \kappa_* := \min \{ r-1, \kappa\}$:
\begin{align}\label{eq:226}
  p_{r,i}^{(j)}(z+wh)=b_{i,\kappa,j}(w)h^{\kappa-j}f^{(\kappa)}(z)+\bigO(h^{\kappa+1-j}).
\end{align}
We use the change of variables $x=z+wh$ to get from \eqref{eq:226} for
$j=1$:
  \begin{gather*}
    \int_{x_{-1/2}}^{x_{1/2}} \bigl(p_{i}^{(1)}(x)\bigr)^2 \mathrm{d}x  
    =h^{2\kappa-1}\mu_{i,1}+\bigO(h^{2\kappa}), \\  
    \mu_{i,1}:= \bigl(f^{(\kappa)}(z)^2 \bigr)\int_{\mz}^{\mo} b_{i,\kappa,1}(w)^2\mathrm{d}w > 0. 
  \end{gather*}
For $1<j\leq  \kappa$  (and, a fortiori, $r>2$, therefore $\kappa=k+1$) we obtain 
  \begin{gather*}
    \int_{x_{-1/2}}^{x_{1/2}} \bigl(p_{r,i}^{(j)}(x) \bigr)^2 \, \mathrm{d}x     
    =\mu_{i,j}h^{2(\kappa-j)+1}+\bigO(h^{2(\kappa-j)+2}), \\
        \mu_{i,j}:= \bigl( f^{(\kappa)}(z) \bigr) ^2\int_{\mz}^{\mo}\bigl(       b_{i,\kappa,j}(w) \bigr)^2 \mathrm{d}w\geq 0. 
  \end{gather*}
For $j > \kappa$  we get 
  \begin{align*}
    \int_{x_{-1/2}}^{x_{1/2}}\bigl(p_{r,i}^{(j)}(x) \bigr)^2 \,\mathrm{d}x  &=
    h
    \int_{x_{-1/2}}^{x_{1/2}}\bigl(  b_{i,j,j}(w)f^{(j)}(z)+\bigO(h)
\bigr)^2 \mathrm{d}w=\bigO(h).
\end{align*}
The proof is complete after substituting  these terms into  \eqref{eq:SI}: 
  \begin{align*}
    I_{i} 
    &=
    \sum_{j=1}^{\kappa_*}h^{2j-1} \big(\mu_{i,j}h^{2(\kappa-j)+1}+\bigO(h^{2(\kappa-j)+2})\big)
    +
    \sum_{j=\kappa_*+1}^{r-1}h^{2j-1}\bigO(h) \\ & =h^{2\kappa}\sum_{j=1}^{\kappa_*}\mu_{i,j}+\bigO(h^{2\kappa+1}),
  \end{align*}
  where we take into account that $\mu_{i,1} + \dots + \mu_{i,\kappa_*} >0$. 
\end{proof}

\section{Design of WENO weights}\label{wd}
 To define our modified scheme (the OWENO scheme), 
   we design weights in such a way that the resulting scheme
 has the order of accuracy $2r-1$, for $r>2$,  corresponding
 to WENO reconstructions of order at least~$5$. We do not consider  the case $r=2$ since 
   severe technical difficulties arise in the accuracy analysis, according
   to the results drawn in Theorem~\ref{is}. This issue is very
   complex to address and will be tackled in full detail in a separate
   paper.
 
 In WENO schemes, the weights $\omega_i$  are defined by a relation of the type
\begin{equation}\label{weights}
  \omega_{i}= \alpha_{i} / (\alpha_{0} + \cdots + \alpha_{r-1}),\quad
  0\leq i\leq r-1, 
\end{equation}
so that  $\omega_0 + \dots + \omega_{r-1}=1$. 
In this section the quantities $\alpha_{0},  \dots ,\alpha_{r-1}$ are given by 
\begin{equation}\label{alpha}
  \alpha_{i}=c_{i}\left(1+\frac{d}{I_{i}^{s_1}+\eps}\right)^{s_2},\quad
  0\leq i\leq r-1, 
\end{equation}
for some $s_1,s_2>0$, $c_i>0$ with $c_0 + \dots + c_{r-1} =1$ and
where $d$~is a function,   to be defined
below, that depends  on $f_{-r+1},\dots,f_{r-1}$. This
approach is related to     Yamaleev and Carpenter
\cite{YamaleevCarpenter2009}.
The ultimate goal is to obtain the order of convergence $2r-1$,
regardless of the presence of neighboring extrema \cite{SINUM2011,ArandigaMartiMulet2014,HenrickAslamPowers2005,YamaleevCarpenter2009},  
and without assuming anything about the small number $\eps>0$ that ensures the
strict positivity of the denominators. In contrast to  other approaches \cite{SINUM2011,ArandigaMartiMulet2014}, 
    our design does {\em not} rely
on a functional relation between $\eps$ and $h$. Although $\eps >0$  is necessary if conditionals are to be avoided (which in turn may be necessary 
 to avoid divisions by zero), 
our arguments will show that $\eps$ can be neglected in the
asymptotical analysis  of the order with respect to $h$.

\subsection{Motivation}\label{subsec:mot}
	 
In the classical WENO order-enhancing argument in case of sufficient
smoothness, for a function with an extremum of order~$k$, the order of the
reconstruction is
\begin{align} \label{ordmax} 
  \text{ord}_{\max}=\min \bigl\{\max \{2r-1, k+1\}, s+\max\{r, k+1\} \bigr\},
\end{align}
where
$\max\{2r-1, k+1\}$, resp.~$\max\{r, k+1\}$, are the orders of the
reconstructions with~$p_{2r-1,r-1}$, resp.~$p_{r,i}$ (see
Lemma~\ref{lemma:1}) and $s\geq 0$ satisfies
$\omega_i=c_i+\bigO(h^s)$. In what follows,  we may assume $k\leq 2r-3$, 
 since otherwise \eqref{ordmax} stipulates
  $\text{ord}_{\max}=2r-1$. 

Yamaleev and Carpenter propose in  \cite{YamaleevCarpenter2009} the
following squared undivided difference of the $2r-1$ consecutive
values $\{f_{-r+1}, \dots, f_{r-1}\}$ to be used in
  \eqref{alpha} as term~$d$: 
\begin{equation}\label{und_dif}
  d := d_1 := \Delta_{2r-2}(f_{-r+1},\dots,f_{r-1}):=
	\Biggl( \, \sum_{j=-r+1}^{r-1}(-1)^{j+r-1}\binom{2r-2}{j+r-1}f_j\Biggr)^2,
      \end{equation}
which has the following asymptotic accuracy properties: 
        \begin{align}\label{eq:pep5} \begin{split} 
  \Delta_{2r-2} \bigl(f_{-r+1},\dots,f_{r-1}\bigr)&=\begin{cases}
  \bigO(h^{4r-4}) & \text{if $f\in C^{2r-2}(z)$,} \\
  \bar{\bigO}(1) & \text{if $f\notin C^{0}(z).$} 
\end{cases}
\end{split}
\end{align} 
Under the smoothness assumption, 
if we  set $d=d_1^{s_1}$ in \eqref{alpha}, then in view of 
$\smash{I_j=\bar\bigO(h^{2k+2})}$ (cf.\  Theorem~\ref{is}) we obtain  
 $\smash{d_1^{s_1}/I_{i}^{s_1}=\bigO(h^{s_1(4r-2k-6)})}$.
The order-enhancing argument in this context requires that 
 $\smash{d_1^{s_1} / I_{i}^{s_1} \to 0}$   as $h \to 0$, 
which is not met if  $k = 2r-3$.  On the other hand,  if
$k\geq 2r-2$, then $\text{ord}_{\max}\geq 2r-1$. So there remains an
order loss gap at $k=2r-3$.  
 We herein close this gap  by proposing   an expression 
$d=D_r^{s_1}$, where the function $D_r = \Delta_{2r-2}(f_{-r+1},\dots,f_{r-1}; \eps)$
is designed such that the second-degree homogeneity property holds
\begin{align}\label{eq:pep2}
  & \Delta_{2r-2}(\alpha f_{-r+1},\dots,\alpha f_{r-1}; 0)  = 
  \alpha^2 \Delta_{2r-2}( f_{-r+1},\dots, f_{r-1}; 0)  \quad \text{for all $\alpha \in \mathbb{R}$,} 
  \end{align} 
  and  that whenever $z_h=z+\bigO(h)$,
  \begin{align} \label{eq:pep1} 
\Delta_{2r-2}\bigl(f_{-r+1},\dots,\alpha f_{r-1};0^+ \bigr)
& =
\begin{cases} \bigO(h^{4r-4}) & \text{if $f\in C^{2r-2}(z)$ and $k < 2r-3$,} \\
  \bigO(h^{4r-3}) & \text{if $f\in C^{2r-2}(z)$ and  $k= 2r-3$,} \\ 
  \bar\bigO(1) & \text{if $f\notin C^0(z),$} 
\end{cases}  
\end{align}
where
  $\Delta_{2r-2}(\cdot;0^+):=\lim_{\varepsilon\to0^+}\Delta_{2r-2}(\cdot;\varepsilon)$.

Clearly, the previous analysis shows that the Yamaleev-Carpenter function~$d_1$  in \eqref{und_dif}
satisfies \eqref{eq:pep2}, but fails to satisfy 
\eqref{eq:pep1} by one order when $k=2r-3$.

\subsection{Novel smoothness indicator}\label{subsec:nsi}
The crucial contribution of this section, and the main novelty of this work,  is  the definition of
a smoothness indicator that satisfies 
\eqref{eq:pep2} and at the same time \eqref{eq:pep1}, namely, behaves
like   $\bigO(h^{4r-3})$, i.e., one
order more than~$d_1$, when $f\in C^{2r-2}(z)$ and~$k= 2r-3$. 
 This new  smoothness indicator is   defined  by 
\begin{align}\label{eq:d2} 
 d_2:= \Delta_{2r-2}(f_{-r+1,h},\dots,f_{r-1,h}): =B_h-4A_hC_h,
\end{align} 
where $A_h$, $B_h$ and $C_h$  are  the coefficients of the
parabola 
\begin{align*} 
P_h^{(2r-4)}(w)=A_hw^2+B_hw+C_h, 
\end{align*} 
which is the $(2r-4)$-th derivative of 
 $P_h(w)=p_h(z+wh)$, where $p_h  \in \Pi_{2r-2}$ is the  reconstruction
polynomial associated with  the data
$f_{-r+1,h},\ldots,f_{r-1,h}$ and $f_{j,h}:=z_h+jh$. Further details on the representation of 
 the derivatives of~$P_h$  can be found in Lemma
\ref{linear}. We state some properties of this new smoothness indicator prior to the definition of 
	the parameter $d$ in \eqref{alpha}.

\begin{proposition}\label{roots_discontinuity}
  Let $n\geq3$. With the
  same notation as in Lemma~\ref{linear}, 
if $f\in C^0(z^{\pm})$ is discontinuous at~$z$, then 
  $\smash{\Delta_n(f(x_{0,h}),\dots,f(x_{n,h}))=\bar{\bigO}(1)}$. 
\end{proposition}

\begin{proof}
 We let
  $f(z^{-})=:f_{\mathrm{L}}\neq
  f_{\mathrm{R}}:=f(z^{+})$, where
    $f(z^{\pm}):=\lim_{y\to z^{\pm}}f(y)$, and define 
  \[i_0:=
  \begin{cases}
    \min\{0\leq i\leq n\mid a_i\leq0\wedge a_{i+1}>0\} &
    \textnormal{if $f(z)=f_{\mathrm{L}}$}, \\
    \min\{0\leq i\leq n\mid a_i<0\wedge a_{i+1}\geq0\} &
    \textnormal{if $f(z)=f_{\mathrm{R}}$,} 
  \end{cases}
   \quad 
    f_i:=\begin{cases}
f_{\mathrm{L}} & \text{if $i\leq i_0$,} \\
f_{\mathrm{R}} & \text{if $i> i_0$.} 
\end{cases}
\] 
 If $p \in  \bar{\Pi}_n$ is   the interpolating polynomial  
 with  $p(z+a_ih)=f_i$,
  $0\leq i\leq n$ and $P(w):=p(z+wh)$, then, by Lemma
  \ref{roots_discont}, $P^{(n-2)}$ has two simple roots,  and therefore  
      $\Delta_n(f_0,\dots,f_n) > 0$. Since 
    $\Delta_n$ is a continuous function (quadratic
      function with respect to their arguments)
  and $\lim_{h\to 0}f(x_{i,h})=f_i$,    
  \begin{align*}
    \lim_{h\to 0}\Delta_n \bigl(f(x_{0,h}),\dots,f(x_{n,h}) \bigr)
    =\Delta_n(f_0,\dots,f_n) > 0,
  \end{align*}
  hence  $\Delta_n(f(x_{0,h}),\dots,f(x_{n,h}))=\bar\bigO(1)$.
\end{proof}

The following result is presented  for  a more general grid of the
form $z_h+a_ih$, where $z_h$  is assumed to satisfy 
$z_h=z+\bigO(h)$. This generalization implies that 
$\Delta_n$ satisfies the desired bounds not only when the critical
point is located in a relative position with respect to the stencil, but
also when the stencil converges
to the critical point (regardless of the relative position with
respect to the critical point) as $h \to 0$.  Namely,
the following result stands for the behaviour of $\Delta_n$
\textit{near} a critical point. This consideration is crucial in the
context of partial differential equations (PDEs), in which the relative position of a critical point
with respect to the stencils selected from the grid is arbitrary.
\begin{proposition}\label{roots_smoothness}
  Let $n\geq3$ and  assume that $f\in C^{n+1}(z)$ satisfies 
  $f^{(n-1)}(z)=f^{(n-2)}(z)=0$,  $f^{(n)}(z)\neq0$. Let
  $z_h\in\R$ such that $z_h-z=\bigO(h)$ and the stencil $x_{i,h}=z_h+a_ih$, $0\leq i\leq n$, $a_0<a_1<\dots<a_n$. Then  there holds 
  $$\Delta_n \bigl(f(x_{0,h}),\dots,f(x_{n,h}) \bigr)=\bigO(h^{2n+1}).$$
\end{proposition}

\begin{proof}
  By Lemma \ref{linear}, there holds
  \begin{align*}
    P_h^{(n-2)}(w)=\sum_{j=0}^2L_{\boldsymbol{a}}^{n-2,j} \bigl(f(x_{0,h}),\ldots,f(x_{n,h}) \bigr)w^j,
  \end{align*}
  where $L_{\boldsymbol{a}}^{n-2,j} = L_{\boldsymbol{a}}^{n-2,j}(f(x_{0,h}),\ldots,f(x_{n,h}))$, $j=0,1,2$,  satisfy 
  $$L_{\boldsymbol{a}}^{n-2,j} \bigl(f(x_{0,h}),\ldots,f(x_{n,h}) \bigr)=\frac{1}{j!}h^{n-2+j}f^{(n-2+j)}(z_h
  )+\bigO(h^{n+1}),\quad j=0,1,2.$$
  Denoting $\delta_h:=z_h-z=\bigO(h)$, 
  $\smash{A_h:=L_{\boldsymbol{a}}^{n-2,2}}$, 
  $\smash{B_h:=L_{\boldsymbol{a}}^{n-2,1}}$, and 
  $\smash{C_h :=L_{\boldsymbol{a}}^{n-2,0}}$, 
    using Taylor
  expansion around~$z$ and  considering  that
  $f^{(n-2)}(z)=f^{(n-1)}(z)=0$, we  obtain 
  \begin{align*}
    A_h&=\frac{1}{2}h^nf^{(n)}(z_h)+\bigO(h^{n+1}) 
    =\frac{1}{2}h^nf^{(n)}(z)+\bigO(h^{n+1}),\\
    B_h&=h^{n-1}f^{(n-1)}(z_h) 
   =\delta_hh^{n-1}f^{(n)}(z)+\bigO(h^{n+1}),\\
    C_h&=h^{n-2}f^{(n-2)}(z_h)+\bigO(h^{n+1}) 
 =\frac{1}{2}\delta_h^2h^{n-2}f^{(n)}(z)+\bigO(h^{n+1}).
  \end{align*}
  Therefore, the discriminant of the quadratic equation $\smash{P_h^{(n-2)}(w)=0}$ becomes 
  \begin{align*}
    B_h^2-4A_hC_h & 
    = \delta_h^2 h^{2n-2} f^{(n)}(z)^2 \big(
    (1+\bigO(h^2))^2-(1+\bigO(h^3))(1+\bigO(h))\big) \\
    &= 
     \bigO(h)^2 h^{2n-2} f^{(n)}(z)^2 \bigO(h) = \bigO(h^{2n+1}).
   \end{align*} 
\end{proof}

\begin{theorem}\label{extreme_accuracy}
  Let $n\geq3$,   $z_h\in\R$~such that $z_h-z=\bigO(h)$, and consider the stencil
  $x_{i,h}=z_h+a_ih$, $0\leq i\leq n$, $a_0<a_1<\dots<a_n$.
  Then
  
  \begin{align*} 
  &\Delta_n \bigl(f(x_{0,h}),\dots,f(x_{n,h}) \bigr)=\\ & =\begin{cases}
  \bar{\bigO}(1) & \textnormal{if there exists $h_0>0$ such that  $x_{0,h}<z<x_{n,h}$}  \\
  & \textnormal{for all $0<h<h_0$, and $f$  has a discontinuity at $z$,}\\
  \bigO(h^{2n+1}) & \textnormal{if $f\in C^{n+1}$  with $f^{(l)}(z)=0$ for $1\leq l\leq n-1$ 
  and $f^{(n)}(z)\neq0$.} 
  \end{cases}
  \end{align*} 
\end{theorem}

\begin{proof}
  The result  follows from Propositions~\ref{roots_discontinuity} and~\ref{roots_smoothness}, respectively.
\end{proof}

We can now proceed to the definition of $d$ appearing in \eqref{alpha}, 
	in a way such that the resulting reconstruction also attains optimal order near critical points of order $2r-3$ (and thus 
  of critical points of any order).

Let $p_h  \in \Pi_n$, $n=2r-2$,  be the interpolating polynomial associated to the stencil~$S$ (see \eqref{Sdef}). The
$(n-2)$-th derivative of the polynomial $P_h(w):=p_h(z+wh)$ is a
second-degree polynomial, which can be written as  
$$P_h^{(n-2)}(w)=C_h+B_hw+A_hw^2,$$
where $A_h, B_h, C_h$  are linear functions of $f_{-r+1},\dots
f_{r-1}$. Now, by Theorem \ref{extreme_accuracy} with $n=2r-2$, the
expression \eqref{eq:d2} satisfies
  \begin{align}\label{eq:pep4} \begin{split} 
   \Delta_{2r-2} \bigl(f_{-r+1},\dots,f_{r-1} \bigr) &=
\begin{cases}
  \bigO(h^{4r-3}) & \text{if $f\in C^{2r-2}(z)$, $k=2r-3$,} \\
  \bar{\bigO}(1) & \text{if $f\notin C^{0}(z)$.} 
\end{cases}  
\end{split}
\end{align}
For instance, for a WENO5 reconstruction ($r=3$) from point values these terms can be written as
\begin{align*}
  A_h&=\frac{1}{2}f_{-2}-2f_{-1}+3f_0-2f_1+\frac{1}{2}f_2,\\
  B_h&=-\frac{1}{2}f_{-2}+f_{-1}-f_1+\frac{1}{2}f_2,\\
  C_h&=-\frac{1}{12}f_{-2}+\frac{4}{3}f_{-1}-\frac{5}{2}f_0+\frac{4}{3}f_1-\frac{1}{12}f_2,
\end{align*}
while for reconstructions from cell averages the formula for $C_h$ must be replaced by 
\begin{align*}
  C_h&=-\frac{1}{8}f_{-2}+\frac{3}{2}f_{-1}-\frac{11}{4}f_0+\frac{3}{2}f_1-\frac{1}{8}f_2.
\end{align*} 
Based on   \eqref{eq:pep5} and \eqref{eq:pep4}, we define
the  function  
\begin{equation}\label{eq:pep6}
  D_{r}:=d:=\frac{d_1^{s_1}|d_2|^{s_1}}{d_1^{s_1}+|d_2|^{s_1}+\eps}
\end{equation}
related to the harmonic mean of
$d_1^{s_1}$ and $d_2^{s_1}$. 
 Its limit when $\eps \to 0$,  namely 
\begin{align*} 
  \bar d=\begin{cases} \displaystyle
    \frac{d_1^{s_1}|d_2|^{s_1}}{d_1^{s_1}+|d_2|^{s_1}} & \text{if $d_1d_2\neq 0$,} \\
    0& \text{otherwise,} 
  \end{cases}    
\end{align*}  
satisfies both desired properties, namely \eqref{eq:pep2} and  \eqref{eq:pep1}.

The asymptotics of the weights 
for $\eps\to0$ are analyzed in the
Appendix and are used to obtain the following theorem.

\begin{theorem}\label{th:smooth}
  If $f\in C^{2r-1}(z)$, $r\geq 3$, then 
    \begin{align*} 
    f(x_{1/2})-  q(x_{1/2})=\bigO(h^{2r-1})+\bigO(\eps^{s_2}). \end{align*}  
\end{theorem}
\begin{proof}
  We define 
  $\bar\omega_i:=\lim_{\eps\to 0}
\omega_{i}$ and $\bar q(x):= \bar \omega_0 p_{0} (x) +  \dots + \bar\omega_{r-1} p_{r-1} (x)$.
  The first step in the proof is to use Lemma \ref{lemma:pep1} to get
  for $e(h)=f(x_{1/2})-     \bar{q}(x_{1/2})$
  \begin{align*}
    &f(x_{1/2})-q(x_{1/2}) =
    f(x_{1/2})-
    \bar{q}(x_{1/2}) +\bar{q}(x_{1/2})-q(x_{1/2})
    \\
    & =e(h)+\sum_{i=0}^{r-1} \bigl(\bar\omega_i-\omega_i\bigr)p_{i}(x_{1/2})
    =e(h)+\sum_{i=0}^{r-1}\bigO(\eps^{s_2})\bigO(1)  =e(h)+\bigO(\eps^{s_2}).
\end{align*}
It only remains to prove that
  \begin{align} \label{ehorder} 
  e(h)=\bigO(h^{2r-1}), 
  \end{align} 
  which  will be achieved by analyzing the behavior of $\bar\omega_i $,
  for which we may assume that
  \begin{equation}\label{eq:426}
    \text{there exists $h_0>0$  such that  $I_j(h)\neq 0$ for all $0<h<h_0$ and all~$j$,}  
  \end{equation}    
  since, otherwise, for each $n$ 
there exist  $h_n>0$ 
and $j_n\in \Jr$ with
$$\lim_{n\to \infty}h_n=0, \quad  I_{j_n}(h_n)=0.$$
It follows that $f$ is
constant on the  
points $\{ x_{j, h_n}\}$, $\smash{j=-r+1+j_n,\dots,j_n}$. Therefore there exists  
 $\smash{\{z_n\}_{n \in \mathbb{N}}}$ with $z_n\to z$
with $f'(z_n)=0$.
 A recursive use of Rolle's theorem and continuity yields that
 $f^{(k)}(z)=0$ for any $k=1,\dots,2r-1$, so Lemma \ref{lemma:1}
 yields $e(h)=\bigO(h^{2r-1})$.

 We may assume that the order $k$ of the critical point  $z$, satisfies
  $ k< 2r-2$,  since, otherwise, if $k\geq 2r-2$,  then
  Lemma~\ref{lemma:1} would yield that 
  $e(h)=\bigO(h^{k+1})=\bigO(h^{2r-1})$ as 
   in~\eqref{ehorder}. 
  Under  this assumption and \eqref{eq:426},  from   
  \eqref{eq:425} we obtain 
  \begin{equation}\label{eq:1}
    \bar\omega_{i} =  c_i
 \Biggl( \sum_{j=0}^{r-1} c_j 
 \left( \frac{\beta_j}{ \beta_i} \right)^{s_2}
 \Biggr)^{-1}, \quad \beta_i = 1+\bar d / I_{i}^{s_1}. 
\end{equation}

 Theorem \ref{is} yields   $\smash{I_j=\bar\bigO(h^{2(k+1)})}$.
By \eqref{eq:pep6}, \eqref{eq:pep4} and  \eqref{eq:pep5} (in that order),  
we deduce that $d=\bigO(h^{s_1\nu})$, where $\nu = 4r-4$ if $k<2r-3$ 
 and $\nu= 
    4r-3$ if $k=2r-3$.  
    We analyze \eqref{eq:1} with these estimates:
\begin{align*}
    \left|\frac{\beta_j}{\beta_i}-1\right|&=\frac{\bar d}{1+\bar d/I_i^{s_1}}
    \frac{|I_i^{s_1}-I_j^{s_1}|}{I_i^{s_1}I_j^{s_1}}
\leq
    \frac{\bar d(I_i^{s_1}+I_j^{s_1})}{I_i^{s_1}I_j^{s_1}}
    =
    \frac{O(h^{\nu})\bigO(h^{(2(k+1))s_1})}{\bar{\bigO}(h^{4s_1(k+1)})},
    \end{align*} 
 which means that 
 \begin{align}    
    \label{eq:2}
   \beta_j / \beta_i     
    &=1+\bigO(h^{\zeta}),\quad \zeta:=2s_1(\nu-k-1). 
  \end{align}
  It follows from \eqref{eq:1} that
  \begin{align} \label{omegaiest} 
    \bar\omega_{i} =
c_i \Biggl( \sum_{j=0}^{r-1} c_j 
\left(1+\bigO(h^\zeta)\right)^{s_2}
\Biggr)^{-1} =
c_i
 \Biggl( \sum_{j=0}^{r-1} c_j 
 \bigl(1+\bigO(h^\zeta) \bigr)
 \Biggr)^{-1}=
c_i+\bigO(h^\zeta).
\end{align}
  Using that 
   $\bar\omega_0 + \dots + \bar\omega_{r-1} = c_0 + \dots + c_{r-1}$,
$f(z+h/2)-p_{2r-1,r-1}(z+h/2)=\bigO(h^{2r-1})$, and \eqref{eihest},
 we obtain from  \eqref{omegaiest} 
\begin{align*}
  e(h)= 
  \sum_{i=0}^{r-1}\bar\omega_i e_i (h) 
  &
  =
  \sum_{i=0}^{r-1} \bigl(c_i+\bigO(h^{\zeta}) \bigr)
   \bigl(f(z+ h/2) -p_i (z+h/2) \bigr) 
  \\
  &=
  \sum_{i=0}^{r-1}c_i \bigl(f (z+h/2) -p_i (z+h/2) \bigr)
+
\sum_{i=0}^{r-1}\bigO(h^{\zeta})\bigO \bigl(h^{\max \{r, k+1\}} \bigr)
\\
&=
f(z+h/2)-p(z+h/2)
+
\bigO\bigl(h^{\zeta+\max\{r, k+1 \}}\bigr)
\\
&=
\bigO(h^{2r-1})
+
\bigO \bigl(h^{\zeta+\max \{ r, k+1 \}} \bigr)
=\bigO \bigl(h^{\min \{2r-1, \zeta+\max \{r, k+1\} \}} \bigr).
\end{align*}
Utilizing the definition of~$\zeta$ in \eqref{eq:2}, one can easily verify that 
  $\zeta+\max \{ r, k+1 \} \geq 2r-1$ for all  
 $ k\leq  2r-3$ and $s_1\geq 1$. 
\end{proof}

\begin{remark}
  All these precautions on the possibility of having smoothness
  indicators that vanish asymptotically are not void, since the function
  \begin{align*}
    f(x)=\begin{cases}
      \mathrm{e}^{-1/x^2}& \text{\em for $x>0$,} \\
      0 & \text{\em for $x\leq 0$} 
    \end{cases}
  \end{align*}
  satisfies $f \in C^{\infty}(\mathbb R)$ and
  $f^{(n)}(0)=0$ for all
  $n\in \mathbb{N}$, therefore, for
  $x=0$, it follows  that $I_0(h)=0$  for all~$h>0$.
\end{remark}

\begin{theorem}\label{th:discontinuous}
  If $f$ has a discontinuity at $z$ and is $r$~times continuously
  differentiable in $(z-\delta_0,z)\cup(z,z+\delta_0)$ for some
  $\delta_0>0$ and   is $r$~times continuously differentiable  either at
  $z^-$ or at $z^+$, then
  \begin{align*}
    f(x_{1/2})-q(x_{1/2})=\bigO \bigl(h^{\min \{ r, 2s_1s_2 \}} \bigr)+\bigO(\eps^{s_2}).
  \end{align*} 
\end{theorem}
\begin{proof}
We
  use the same notation and  assume  that
  $\eps=0$ as in the proof of Theorem~\ref{th:smooth} 
 and aim to prove that $\smash{e(h)=\bigO(h^{\min\{r, 2s_1s_2\}})}$. 
 We define the index set 
 \begin{align*}  
 J_r:=\bigl\{ 0\leq j\leq r-1: f|_{[x_{-r+1+j}, x_{j}]}\in C^r  \bigr\}.
 \end{align*} 
  By the assumption on the lateral smoothness of~$f$ at~$z$, since  $z\in [x_{-r+1,j}, x_{j}]$ if and only if
  $-r+1+i\leq (z-z_h)/h
  \leq i$
  and $(z-z_h)/h\in (-1, 1)$, it follows that
  \begin{equation}
    \begin{cases}
      0\in J_r &
      \text{
  if $(z-z_h)/h\in (0, 1)$ or $z=z_h=x_0$ and
  $f\in C^r(z^-)$,}\\
r-1\in J_r &
\text{
  if $(z-z_h)/h\in (-1, 0)$ or $z=z_h=x_0$ and
  $f\in C^r(z^+)$,}
\end{cases}
\end{equation}
hence  $J_r\neq\varnothing$.

 The main difference with respect to Theorem \ref{th:smooth}   is that 
      $\smash{I_j=\bigO(h^{m_j})}$, where 
      $m_j =0$ if $j\notin J_r$ and $m_j=       2(k+1)$ if $j\in J_r$ and     $d=\bar\bigO(1)$, 
which immediately yields
  \begin{align}
  \notag
  \frac{\beta_j}{\beta_i}&=\frac{1+d/I_j^{s_1}}{1+d/I_i^{s_1}} = \bar\bigO \bigl(h^{(m_i-m_j)s_1} \bigr).
  \end{align}
  Therefore, for $i\notin J_r$, \eqref{eq:1} reads
  \begin{align*}
    \bar\omega_i &=c_i \Biggl( \sum_{j\in J_r}c_j
      \left(\frac{\beta_j}{\beta_i}\right)^{s_2}+\sum_{j\not\in J_r}
      c_j \left(\frac{\beta_j}{\beta_i}\right)^{s_2} \Biggr)^{-1}\\ &  =
    c_i \Biggl( \sum_{j\in J_r}c_j
      \bigl(\bar\bigO(h^{-2(k+1)s_1})\bigr)^{s_2}+\sum_{j\not\in J_r}
      c_j \left(\bar\bigO(1)\right)^{s_2} \Biggr)^{-1} \\
    &=
    \frac{c_i}{\bar\bigO(h^{-2(k+1)s_1s_2})+\bar\bigO(1)}=
    \frac{c_i}{\bar\bigO(h^{-2(k+1)s_1s_2})}=
    \bigO \bigl(h^{2(k+1)s_1s_2} \bigr)=
    \bigO \bigl(h^{2s_1s_2} \bigr) 
  \end{align*}
since $k\geq 0$.
  Since  $\bar\omega_i\leq 1$, $ e_i(h)= \bigO(1)$ if $i\notin J_r$ and 
     $e_i(h)= \bigO(h^r)$ if $i\in J_r$, 
  we deduce
  \begin{align*}
  e(h)&= 
  \sum_{i=0}^{r-1}\bar\omega_i e_i (h) 
  =
  \sum_{i\notin J_r}\bigO \bigl( h^{2s_1s_2} \bigr) \bigO(1) +
  \sum_{i\in J_r}\bigO(1) \bigO(h^r) 
  =
  \bigO \bigl(h^{\min \{r,2s_1s_2\} } \bigr).
  \end{align*}
\end{proof}

\begin{remark}\label{s1s2}
  As a consequence of Theorem \ref{th:discontinuous}, we may 
  take $2s_1s_2\geq r$ to get the suboptimal $r$-th order at
  discontinuities.
\end{remark}

\subsection{Summary of the algorithm}\label{sec:alg1}
 For the ease of reference we summarize here 
 the  new OWENO reconstruction for a local stencil. 
 
\smallskip  

\noindent Input: $\{f_{-r+1},\ldots,f_{r-1}\}$ and $\eps>0$.
\begin{enumerate}
\item Compute $p_{i}$, $0\leq i\leq r-1$, the corresponding
  reconstruction polynomials of degree $r-1$ at
  $x=x_{1/2}$. See \cite[Proposition 1]{SINUM2011}
    for further details about their explicit expression.
\item Compute the Jiang-Shu smoothness indicators \eqref{eq:SI}. See
  \cite[Proposition 5]{SINUM2011} for further details about the
  explicit computation procedure to obtain their expression.
\item Compute ${d}$ from \eqref{eq:pep6}  
for $d_1 := \Delta_{2r-2}$ as given by \eqref{und_dif}, and
$d_2:=\Delta_r$ as given in \eqref{eq:d2}.
\item Compute the terms $\alpha_i$ from \eqref{alpha}, where $d$ is given by \eqref{eq:pep6}, with $c_{i}$ the ideal linear weights, for some $s_1,s_2$ chosen by the user such that $s_1\geq1$ and $s_2\geq r/(2s_1)$.
\item Generate the WENO weights $\omega_0, \dots, \omega_{r-1}$ from \eqref{weights}. 
\item Obtain the OWENO reconstruction at $x_{1/2}$:
  \begin{align*}  
  q_r(x_{1/2})= \omega_{0}p_{0}(x_{1/2}) + \dots + \omega_{r-1}p_{r-1}(x_{1/2}).
  \end{align*}
\end{enumerate}
Output: $q_r(x_{1/2})$.

\begin{remark}\label{rem:even}
  Since  it is not guaranteed  that $d_2\geq0$, we included its
  absolute value $|d_2|$ in Equation
  \eqref{eq:pep6}. If one wants to avoid using an absolute value (and thus a Boolean condition in a WENO scheme), one has simply to chose an even $s_1$ satisfying the bounds in Remark \ref{s1s2}.
\end{remark}

\section{Numerical experiments}\label{sec:num_exp}

In this section, the chosen exponents are $s_1=2\lceil r/4\rceil$
(taking into account Remark \ref{rem:even}), and $s_2=1$. The reason
for this choice is  that the choice of $\varepsilon$ in \eqref{alpha}
is related to the exponent $s_2$, since one should take
$\smash{\varepsilon\gtrsim  \varepsilon_0^{1/{s_2}}}$, with
$\varepsilon_0$ the lowest positive number of the working precision,
in order to avoid arithmetic underflow/overflow. Moreover, although
unnecessary according to the accuracy requirements in case of
smoothness, the greater the parameter $s_1$ is, the closer are
simultaneously the weights to the ideal weights in case of smoothness
and to zero in case of discontinuity. 
\subsection{Algebraic test cases}\label{subsec:alg_prob}
We start our numerical tests with several numerical experiments 
devoted to  emphasize the accuracy properties analyzed theoretically
beforehand. We will perform tests involving JS-WENO (with the  weight
design by  Jiang and Shu \cite{JiangShu96}), WENO-Z \cite{CastroCostaDon2011}, YC-WENO \cite{YamaleevCarpenter2009} (with the improved version of the Yamaleev-Carpenter weight design 
 \cite{ArandigaMartiMulet2014}; 
   and OWENO (with our design) reconstructions of order $2r-1$, with $2\leq r\leq 5$. All tests are performed  with  reconstructions both from cell average values to pointwise values and from pointwise values to pointwise values.

We perform these experiments by using the multiple-precision library MPFR \cite{MPFR} through its C++ wrapper \cite{Holoborodko}, using a precision of $3322$ bits ($\approx1000$ digits) and taking  $\smash{\varepsilon=10^{-10^6}}$ in all cases.
\subsubsection*{Example 1: Smooth problem}

\begin{table}[t]
  \setlength\tabcolsep{3pt}
  \footnotesize 
  \begin{center} 
  \begin{tabular}{|c|c|c|c|c|c|c|c|c|}
    \hline
    $k$ & JS-WENO & WENO-Z & YC-WENO & OWENO & JS-WENO & WENO-Z & YC-WENO & OWENO \\
    \hline
    & \multicolumn{4}{|c|}{Order 5 (from point values)} & \multicolumn{4}{|c|}{Order 5 (from cell averages)} \\
    \hline
    0 & 4.9915 & 5.0022 & 4.9983 & 4.9983 & 4.9909 & 5.0018 & 4.9983 & 4.9983 \\
    1 & 3.9742 & 5.0161 & 4.9980 & 4.9980 & 3.9802 & 5.0203 & 4.9981 & 4.9980 \\
    2 & 3.0198 & 2.9777 & 5.0331 & 5.0324 & 3.0348 & 2.9749 & 5.0324 & 5.0317 \\
    3 & {\bf 3.9946} & {\bf 3.9945} & {\bf 3.9945} & {\bf 5.0056} & {\bf 3.9928} & {\bf 3.9927} & {\bf 3.9928} & {\bf 5.0035} \\
    \hline
     & \multicolumn{4}{|c|}{Order 7 (from point values)} & \multicolumn{4}{|c|}{Order 7 (from cell averages)} \\
    \hline
    0 & 6.9902 & 6.9982 & 6.9984 & 6.9984 & 6.9899 & 6.9982 & 6.9984 & 6.9984 \\
    1 & 5.9743 & 7.0023 & 6.9981 & 6.9981 & 5.9699 & 7.0012 & 6.9981 & 6.9981 \\
    2 & 5.0494 & 7.0424 & 7.0002 & 7.0000 & 5.0432 & 7.0363 & 7.0001 & 6.9998 \\
    3 & 4.0005 & 4.0005 & 7.0627 & 7.0548 & 4.0001 & 4.0001 & 7.0600 & 7.0482 \\
    4 & 5.0747 & 5.0747 & 7.0040 & 7.0040 & 5.0655 & 5.0655 & 7.0108 & 7.0108 \\
    5 & {\bf 6.0008} & {\bf 6.0008} & {\bf 6.0008} & {\bf 6.9907} & {\bf 6.0011} & {\bf 6.0011} & {\bf 6.0011} & {\bf 6.9980} \\
    \hline
 & \multicolumn{4}{|c|}{Order 9 (from point values)} & \multicolumn{4}{|c|}{Order 9 (from cell averages)} \\
    \hline
    0 & 8.9831 & 8.9984 & 8.9984 & 8.9984 & 8.9829 & 8.9985 & 8.9985 & 8.9985 \\
    1 & 8.0225 & 8.9983 & 8.9983 & 8.9983 & 8.0226 & 8.9983 & 8.9983 & 8.9983 \\
    2 & 7.0368 & 9.0879 & 8.9981 & 8.9981 & 7.0229 & 9.0782 & 8.9981 & 8.9981 \\
    3 & 6.0712 & 9.0245 & 8.9978 & 8.9978 & 6.0625 & 9.0159 & 8.9978 & 8.9979 \\
    4 & 5.0133 & 5.0133 & 9.0628 & 8.9976 & 5.0072 & 5.0072 & 9.0625 & 8.9976 \\
    5 & 5.9855 & 5.9855 & 9.0325 & 9.0185 & 5.9815 & 5.9815 & 9.0283 & 9.0082 \\
    6 & 7.0409 & 7.0409 & 9.0121 & 9.0121 & 7.0746 & 7.0746 & 9.0143 & 9.0143 \\ 
    7 & {\bf 7.9898} & {\bf 7.9898} & {\bf 7.9898} & {\bf 8.9541} & {\bf 7.9880} & {\bf 7.9880} & {\bf 7.9880} & {\bf 8.9872} \\
    \hline
  \end{tabular}
  \end{center} 
  \smallskip 
  \caption{Example 1 (smooth problem): Fifth-order, seventh-order, and ninth-order reconstructions. The cases in which both JS-WENO and YC-WENO methods lose accuracy (critical point of order $2r-3$) have been highlighted in bold text, in which it can be observed that the  OWENO method keeps the optimal accuracy.}
  \label{smooth1}
\end{table}

Let us consider the family of functions $f_k:\mathbb{R}\to\mathbb{R}$, $k\in\mathbb{N}$, given by
$\smash{f_k(x)=x^{k+1}\mathrm{e}^x}$. 
The function $f_k$ has a critical point at $x=0$ of order $k$. Results involving the different values of $r$ and $k$ considered ($0\leq k\leq 2r-3$) are shown for $3\leq r\leq5$ in 
 Table~\ref{smooth1} for the case of JS-WENO, YC-WENO and OWENO reconstructions. The error is given by $E_{k,n}=|P_N(0)-f_k(0)|$, with $P$ the corresponding reconstruction at $x_{1/2}=0$, with the grid $x_i=(i-1/2)h$, $-r+1\leq i\leq r-1$, with $h=1/N$ for $N\in\N$, when pointwise values 
 \eqref{fpointvalues} (with $f=f_k$) 
 are taken,  and  pointwise values are reconstructed from  pointwise values. Table~\ref{smooth1} 
  also presents the results for the  same setup when cell average values \eqref{fcellaverages}  (with $f=f_k$)  are taken instead
   and pointwise values are reconstructed from cell averages. In all  cases Table~\ref{smooth1}  shows the corresponding average reconstruction orders
   \begin{align*} 
   O_k=\frac{1}{80}\sum_{j=1}^{80}o_{k,j}, \quad \text{where} \quad o_{k,j}=\log_2 \biggl( \frac{E_{k,N_{j-1}}}{E_{k,N_j}}\biggr), 
  \quad 
   N_j=5\cdot2^j, \quad 0\leq j\leq80.
\end{align*} 

As we can see, the JS-WENO loses accuracy near critical points, presenting the order $r+|k-r+1|$, with $k$ the order of the critical point; also, WENO-Z presents the optimal $(2r-1)$-th order for $k<r-1$ and drops to order $k+1$ if $k\geq r-1$, whereas the YC-WENO reconstruction loses accuracy in the corner case $k=2r-3$, as suggested in our theoretical considerations. In contrast, the OWENO reconstructions attain the optimal accuracy in all cases. This confirms that in practice  the OWENO reconstruction  is indeed able to overcome the loss of accuracy in all cases, including those in which YC-WENO-type reconstructions fail to attain the optimal accuracy. 

\subsubsection*{Example 2: Discontinuous problem}

\begin{table}[t]
 \begin{center} 
 {\setlength\tabcolsep{3pt}\footnotesize 
  \begin{tabular}{|c|c|c|c|c|c|c|c|c|}
    \hline
     $\theta$ & JS-WENO & WENO-Z & YC-WENO & OWENO & JS-WENO & WENO-Z & YC-WENO & OWENO \\
    \hline
    & \multicolumn{4}{|c|}{Order 5 (from point values)} & \multicolumn{4}{|c|}{Order 5 (from cell averages)} \\
    \hline
    -2 & 2.9955 & 2.9952 & 2.9951 & 2.9917 & 2.9955 & 2.9952 & 2.9951 & 2.9929 \\
    -1 & 2.9927 & 2.9912 & 2.9925 & 2.9923 & 2.9935 & 2.9925 & 2.9934 & 2.9933 \\
    0 & 3.0029 & 3.0081 & 3.0045 & 3.0070 & 3.0033 & 3.0088 & 3.0050 & 3.0071 \\
    1 & 3.0271 & 3.0458 & 3.0390 & 3.0517 & 3.0294 & 3.0478 & 3.0411 & 3.0517 \\
    \hline
    & \multicolumn{4}{|c|}{Order 7 (from point values)} & \multicolumn{4}{|c|}{Order 7 (from cell averages)} \\
    \hline
     -3 & 3.9970 & 3.9980 & 4.0035 & 4.0140 & 3.9971 & 3.9982 & 4.0041 & 4.0297 \\
    -2 & 4.0088 & 4.0091 & 4.0089 & 4.0090 & 4.0071 & 4.0074 & 4.0072 & 4.0073 \\
    -1 & 3.9509 & 3.9487 & 3.9493 & 3.9473 & 4.0086 & 3.9479 & 4.0087 & 4.0088 \\
    0 & 4.0086 & 3.9412 & 4.0086 & 4.0086 & 4.0086 & 3.9407 & 4.0086 & 4.0086 \\
    1 & 4.0234 & 4.0234 & 4.0234 & 4.0234 & 4.0234 & 4.0235 & 4.0234 & 4.0235 \\
    2 & 4.0206 & 4.0257 & 4.0368 & 4.0344 & 4.0211 & 4.0261 & 4.0370 & 4.0353 \\
    \hline
& \multicolumn{4}{|c|}{Order 9 (from point values)} & \multicolumn{4}{|c|}{Order 9 (from cell averages)} \\
    \hline
    -4 & 4.9937 & 4.9937 & 4.9937 & 4.9937 & 4.9938 & 4.9938 & 4.9938 & 4.9938 \\
    -3 & 4.9933 & 4.9933 & 4.9933 & 4.9933 & 4.9933 & 4.9933 & 4.9933 & 4.9933 \\
    -2 & 4.9928 & 4.9928 & 4.9928 & 4.9928 & 4.9927 & 4.9927 & 4.9927 & 4.9927 \\
    -1 & 4.9925 & 4.9924 & 4.9925 & 4.9825 & 4.9924 & 4.9923 & 4.9924 & 4.9924 \\
    0 & 4.9886 & 5.0631 & 4.9886 & 4.9886 & 4.9917 & 5.0634 & 4.9917 & 4.9917 \\
    1 & 5.0561 & 5.0561 & 5.0561 & 5.0561 & 5.0561 & 5.0561 & 5.0561 & 5.0561 \\
    2 & 5.0564 & 5.0564 & 5.0564 & 5.0564 & 5.0574 & 5.0574 & 5.0574 & 5.0574 \\
    3 & 5.0129 & 5.0356 & 5.0992 & 5.1073 & 5.0154 & 5.0373 & 5.1006 & 5.1042 \\
    \hline
  \end{tabular}}
  \end{center} 
  \smallskip 
  \caption{Example 2 (discontinuous problem): Fifth-order, seventh-order, and ninth-order reconstructions. The optimal accuracy is kept by all the reconstructions regardless of the location of the discontinuity.}
  \label{discontinuous1}
\end{table}

We next  test the accuracy of the methods with the same parameters as above  for  the function  
\[f(x)=
\begin{cases}
  \mathrm{e}^x & \text{if $x\leq0$,}  \\
  \mathrm{e}^{x+1} & \text{if $x>0$,} 
\end{cases}
\]
where, in order to highlight the behaviour of the OWENO reconstructions  at discontinuities, we change the location of the discontinuity by 
 utilizing   a grid of the form $x_i=(i-1/2+\theta)h$, $-r+1\leq i\leq r-1$, for $-r+2\leq\theta\leq r-1$. Since $x_{1/2}=\theta h$, the error  is now  given by $|P(\theta h)-g(\theta h)|$. The results are shown in Table~\ref{discontinuous1}. 
  Clearly,   the suboptimal $r$-th order accuracy is also attained in all the cases when the data contain a discontinuity.

\subsection{Experiments for conservation laws}\label{subsec:cons_laws}
In this section some numerical experiments involving hyperbolic conservation laws will be considered. For this purpose, we use a local Lax-Friedrichs
 (LLF) type flux splitting \cite{shuosher88} for smooth problems, and
 Donat-Marquina's flux formula \cite{DonatMarquina96}   for problems with weak solutions.
  On the other hand, for the time discretization, the approximate Lax-Wendroff schemes proposed by 
   Zor\'{\i}o  et al.\  \cite{ZorioBaezaMulet17}   matching the spatial order will be considered. In this section 
  we work in all experiments with double precision representation and set $\varepsilon=10^{-100}$. For all schemes we consider  fifth-order accuracy.

\subsubsection*{Example 3: Linear advection equation}

\begin{table}
  \centering
  {\footnotesize  \addtolength{\tabcolsep}{-2.5pt} 
  \begin{tabular}{|c|c|c|c|c|c|c|c|c|}
    \hline   
     & \multicolumn{2}{c|}{$\|\cdot\|_1$} & \multicolumn{2}{c|}{$\|\cdot\|_{\infty}$} & \multicolumn{2}{c|}{$\|\cdot\|_1$} & \multicolumn{2}{c|}{$\|\cdot\|_{\infty}$} \\  \hline 
   $N$ & Error & rate & Error & rate & Error & rate & Error & rate  \\
    \hline  
    & \multicolumn{4}{c|}{JS-WENO5} & \multicolumn{4}{c|}{WENO-Z5}  \\
    \hline  
    10 & 8.44e-03 & --- & 1.28e-02 & --- & 1.22e-03 & --- & 1.99e-03 & ---  \\
    20 & 3.59e-04 & 4.56 & 6.93e-04 & 4.20 & 3.27e-05 & 5.21 & 5.25e-05 & 5.24  \\
    40 & 1.09e-05 & 5.04 & 2.37e-05 & 4.87 & 1.01e-06 & 5.01 & 1.99e-03 & 5.04  \\
    80 & 3.29e-07 & 5.05 & 7.00e-07 & 5.08 & 3.15e-08 & 5.01 & 4.94e-08 & 5.01  \\
    160 & 1.02e-08 & 5.01 & 2.21e-08 & 4.98 & 9.79e-10 & 5.01 & 1.54e-09 & 5.01  \\
    320 & 3.19e-10 & 5.00 & 6.65e-10 & 5.06 & 3.05e-11 & 5.00 & 4.79e-11 & 5.00  \\
    640 & 9.96e-12 & 5.00 & 2.02e-11 & 5.04 & 9.52e-13 & 5.00 & 1.50e-12 & 5.00  \\
    \hline
     & \multicolumn{4}{c|}{YC-WENO5} & \multicolumn{4}{c|}{OWENO5} \\ \hline 
    10 & 1.02e-03 & --- & 1.55e-03 & --- & 9.52e-04 & --- & 1.45e-03 & --- \\
    20 &  3.27e-05 & 4.96 & 5.16e-05 & 4.91 & 2.95e-05 & 5.01 & 4.65e-05 & 4.96 \\
    40 &  1.01e-06 & 5.01 & 1.60e-06 & 5.01 & 9.03e-07 & 5.03 & 1.42e-06 & 5.03 \\
    80 &  3.15e-08 & 5.01 & 4.94e-08 & 5.01 & 2.78e-08 & 5.02 & 4.37e-08 & 5.02 \\
    160 &  9.79e-10 & 5.01 & 1.54e-09 & 5.01 & 8.63e-10 & 5.01 & 1.36e-09 & 5.01 \\
    320 &  3.05e-11 & 5.00 & 4.79e-11 & 5.00 & 2.68e-11 & 5.01 & 4.22e-11 & 5.01 \\
    640 &  9.52e-13 & 5.00 & 1.50e-12 & 5.00 & 8.37e-13 & 5.00 & 1.32e-12 & 5.00 \\
    \hline
  \end{tabular}} 
  \smallskip 
  \caption{Example~3 (linear advection equation, solution at $\smash{T=1}$): fifth-order schemes.}
  \label{advlin_o5}
\end{table}

We consider the linear advection equation with the following domain, boundary condition and initial condition:
\begin{align*} 
& u_t+f(u)_x=0,\quad\Omega=(-1,1),\quad u(-1, t)=u(1, t), \\
& f(u)=u, \quad u_0(x)=0.25+0.5\sin(\pi x), \end{align*} 
whose exact solution is $u(x,t)=0.25+0.5\sin(\pi(x-t))$.
We run several simulations with final time $T=1$, for resolutions $h=2/N$,  $N\in \mathbb{N}$, using the classical JS-WENO, WENO-Z and YC-WENO schemes and the OWENO schemes, and compare them for the case of fifth-order accuracy, both with the $L^1$ and $L^{\infty}$ errors. Since the characteristics point to the right, we use left-biased reconstructions. The results are shown in Table~\ref{advlin_o5} for the fifth-order schemes. All schemes keep  fifth-order accuracy. The results of the OWENO schemes are almost identical to those of the YC-WENO scheme.

\subsubsection*{Examples 4 and 5: Burgers equation}

\begin{table}[t] 
  \centering
  {\footnotesize  
   \addtolength{\tabcolsep}{-2.5pt} 
  \begin{tabular}{|c|c|c|c|c|c|c|c|c|}
    \hline   
     & \multicolumn{2}{c|}{$\|\cdot\|_1$} & \multicolumn{2}{c|}{$\|\cdot\|_{\infty}$} & \multicolumn{2}{c|}{$\|\cdot\|_1$} & \multicolumn{2}{c|}{$\|\cdot\|_{\infty}$} \\  \hline 
   $N$ & Error & rate & Error & rate & Error & rate & Error & rate  \\
    \hline  
    & \multicolumn{4}{c|}{JS-WENO5} & \multicolumn{4}{c|}{WENO-Z5}  \\
    \hline  
    40 & 6.28e-05 & --- & 2.73e-04 & --- & 7.99e-05 & --- & 2.44e-04 & ---  \\
    80 & 3.14e-06 & 4.32 & 4.26e-05 & 2.68 & 6.08e-06 & 3.72 & 3.64e-05 & 2.75  \\
    160 & 1.55e-07 & 4.35 & 2.87e-06 & 3.89 & 4.05e-07 & 3.94 & 4.76e-06 & 2.94  \\
    320 & 9.44e-09 & 4.03 & 2.75e-07 & 3.38 & 2.63e-08 & 3.94 & 5.86e-07 & 3.02  \\
    640 & 5.38e-10 & 4.13 & 3.29e-08 & 3.06 & 1.66e-09 & 3.98 & 6.99e-08 & 3.07  \\
    1280 & 3.46e-11 & 3.96 & 3.58e-09 & 3.20 & 1.03e-10 & 4.01 & 8.22e-09 & 3.09  \\
    2560 & 2.10e-12 & 4.04 & 4.80e-10 & 2.90 & 6.37e-12 & 4.02 & 9.60e-10 & 3.10  \\
    \hline
    & \multicolumn{4}{c|}{YC-WENO5} & \multicolumn{4}{c|}{OWENO5} \\ \hline 
    40 &  2.55e-05 & --- & 2.62e-04 & --- & 2.49e-05 & --- & 2.62e-04 & --- \\
    80 &  8.46e-07 & 4.91 & 1.04e-05 & 4.65 & 8.46e-07 & 4.88 & 1.04e-05 & 4.65 \\
    160 &  2.62e-08 & 5.01 & 3.27e-07 & 4.99 & 2.62e-08 & 5.01 & 3.27e-07 & 4.99 \\
    320 &  7.97e-10 & 5.04 & 1.02e-08 & 5.00 & 7.97e-10 & 5.04 & 1.02e-08 & 5.00 \\
    640 &  2.45e-11 & 5.02 & 3.14e-10 & 5.02 & 2.45e-11 & 5.02 & 3.14e-10 & 5.02 \\
    1280 &  7.59e-13 & 5.01 & 9.71e-12 & 5.02 & 7.59e-13 & 5.01 & 9.71e-12 & 5.02 \\
    2560 &  2.34e-14 & 5.02 & 3.03e-13 & 5.00 & 2.34e-14 & 5.02 & 3.03e-13 & 5.00 \\
    \hline
  \end{tabular}} 
  \smallskip 
  \caption{Example 4 (Burgers equation, smooth solution at $T=0.3$): fifth-order schemes.}
  \label{burgers_o5}
\end{table}

\begin{figure}[t] 
  \centering
     \includegraphics{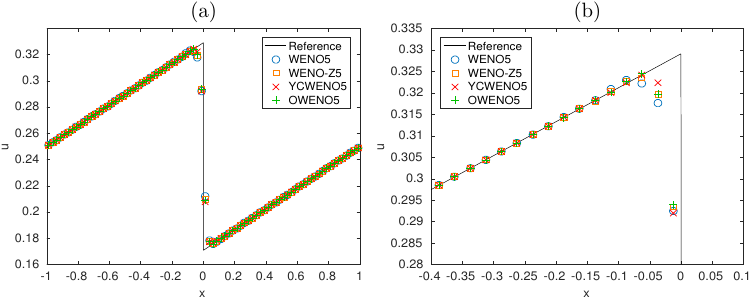} 
  \caption{Example 5 (Burgers equation, discontinuous solution at $T=12$): 
   fifth-order schemes.}
  \label{burgersdisc_o5}
\end{figure}


We now consider the inviscid Burgers equation  along with the following  boundary  and initial conditions:
\begin{align} \label{burgersprob} \begin{split} 
&  u_t+f(u)_x=0,\quad\Omega=(-1,1),\quad u(-1, t)=u(1, t), \\
& f(u)=0.5 u^2, \quad 
u_0(x)=0.25+0.5\sin(\pi x).
 \end{split} 
\end{align}
In this case, $f(u_0(x))$ has a first-order critical point at $x=-1/2$ and $x=1/2$. 
In  Example~4, we consider  the solution of \eqref{burgersprob} at  $T=0.3$, 
 when it remains smooth, while in Example~5 we set $T=12$, when the solution of \eqref{burgersprob} 
  has become discontinuous.   In Example~4 we  run  simulations  for  several  resolutions, with an LLF   
    flux splitting, and 
  display the behaviour of the fifth-order schemes in 
Table~\ref{burgers_o5}. The exact solution is computed through a characteristic line method together with the Newton method, setting as 
 tolerance double-precision machine accuracy.
A loss of the order of  accuracy  is observed for the JS-WENO and WENO-Z schemes. In contrast, the order of  accuracy  of the YC-WENO and all the OWENO schemes is optimal.

In Example~5 we  run the simulation instead until $T=12$. At $t=1$,
the wave breaks and a shock is generated. Therefore, in this case we
use the Donat-Marquina flux-splitting algorithm \cite{DonatMarquina96}. The results are shown in Figure~\ref{burgersdisc_o5} with a resolution of $N=80$ points, and are compared against a reference solution computed with $N=16000$. This ranking of resolution is also consistent with the results for the smooth case.

\subsubsection*{Example 6: Customized equation with a third-order zero}

\begin{table}[t] 
  \centering
  {\footnotesize \addtolength{\tabcolsep}{-2.5pt} 
  \begin{tabular}{|c|c|c|c|c|c|c|c|c|}
    \hline   
     & \multicolumn{2}{c|}{$\|\cdot\|_1$} & \multicolumn{2}{c|}{$\|\cdot\|_{\infty}$} & \multicolumn{2}{c|}{$\|\cdot\|_1$} & \multicolumn{2}{c|}{$\|\cdot\|_{\infty}$} \\  \hline 
   $N$ & Error & rate & Error & rate & Error & rate & Error & rate  \\
    \hline  
    & \multicolumn{4}{c|}{JS-WENO5} & \multicolumn{4}{c|}{WENO-Z5}  \\
    \hline  
    40 & 7.96e-05 & --- & 5.17e-04 & --- & 6.94e-05 & --- & 5.14e-04 & ---  \\
    80 & 4.67e-06 & 4.09 & 7.31e-05 & 2.82 & 3.81e-06 & 4.19 & 7.29e-05 & 2.82  \\
    160 & 2.70e-07 & 4.11 & 9.73e-06 & 2.91 & 2.18e-07 & 4.13 & 9.70e-06 & 2.91  \\
    320 & 1.60e-08 & 4.08 & 1.25e-06 & 2.96 & 1.31e-08 & 4.05 & 1.25e-06 & 2.95  \\
    640 & 9.70e-10 & 4.04 & 1.59e-07 & 2.98 & 8.06e-10 & 4.03 & 1.59e-07 & 2.98  \\
    1280 & 5.95e-11 & 4.03 & 2.01e-08 & 2.99 & 4.99e-11 & 4.01 & 2.00e-08 & 2.99  \\
    2560 & 3.68e-12 & 4.02 & 2.52e-09 & 2.99 & 3.10e-12 & 4.01 & 2.51e-09 & 2.99  \\
    \hline
     & \multicolumn{4}{c|}{YC-WENO5} & \multicolumn{4}{c|}{OWENO5} \\ \hline
    40 &  4.97e-05 & --- & 3.15e-04 & --- & 2.93e-05 & --- & 2.01e-04 & --- \\
    80 &  2.88e-06 & 4.11 & 5.58e-05 & 2.50 & 1.01e-06 & 4.86 & 9.83e-06 & 4.35 \\
    160 &  1.69e-07 & 4.09 & 7.98e-06 & 2.81 & 3.05e-08 & 5.05 & 3.42e-07 & 4.85 \\
    320 &  1.01e-08 & 4.06 & 1.06e-06 & 2.91 & 8.82e-10 & 5.11 & 1.22e-08 & 4.81 \\
    640 &  6.16e-10 & 4.03 & 1.36e-07 & 2.96 & 2.61e-11 & 5.08 & 3.95e-10 & 4.95 \\
    1280 &  3.81e-11 & 4.02 & 1.72e-08 & 2.98 & 7.91e-13 & 5.04 & 1.25e-11 & 4.99 \\
    2560 &  2.36e-12 & 4.01 & 2.17e-09 & 2.99 & 2.51e-14 & 4.98 & 3.90e-13 & 5.00 \\
    \hline
    
  \end{tabular}} 
  \smallskip 
  \caption{Example 6 (customized equation, smooth solution at $T=0.3$): fifth-order schemes.}
  \label{custom_o5}
\end{table}

We now consider the following initial-boundary value problem 
 for a customized equation:
\begin{align*} 
 & u_t+f(u)_x=0,\quad\Omega=(-1,1),\quad u(-1, t)=u(1, t), \\
& f(u)= 0.5 u^2 +0.25u, \quad 
 u_0(x)=0.25+0.5\sin(\pi x).
  \end{align*}  
In this case, $f(u_0(x))$ has a third-order critical point at $x=-1/2$ and a first-order critical point at $x=1/2$.
We now compare the behaviour of the three schemes with the same setup as in  Example~4, by running a simulation until time $T=0.3$, at which the solution is smooth. For the computation of the exact solution, we once again use the method of characteristic lines, with a Newton method matching the machine accuracy for the double precision. Since in this case the characteristics point always to the right, we use a left-biased upwind scheme. The results are shown in Table~\ref{custom_o5} for the fifth-order schemes.
Clearly,  the optimal order of accuracy is lost for both the JS-WENO, WENO-Z and YC-WENO schemes. In contrast, the fifth-order accuracy is solidly kept by the OWENO schemes. This is another confirmation, this time in the context of conservation laws, in which the OWENO are capable to handle the case $k=2r-3$, unlike the previously existing WENO schemes. 

\subsubsection*{Example 7: Shu-Osher problem}

\begin{figure}[t]
  \centering
  \includegraphics{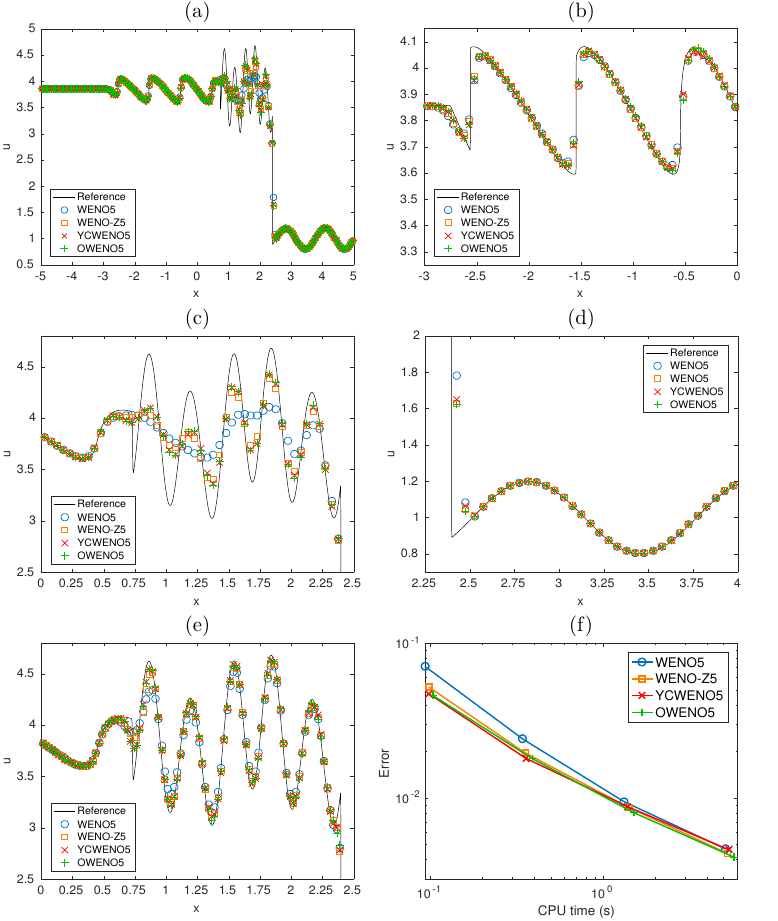} 
  \caption{Example~7 (Euler equations, Shu-Osher problem): numerical solutions at $T=1.8$ by fifth-order schemes:  
   (a)  simulated density for spatial discretization $N=200$, (b--d) enlarged views, (e)   simulated density  for   $N=400$, (f) efficiency plot.}
  \label{shuosher_o5_n=200}
\end{figure}


The 1D Euler equations for gas dynamics are given by 
 $\boldsymbol{u} = ( \rho, \rho v, E)^{\mathrm{T}}$ and $\boldsymbol{f} (\boldsymbol{u}) = 
  \boldsymbol{f}^1  (\boldsymbol{u}) = (\rho v, 
        p+\rho v^2, 
        v(E+p))^{\mathrm{T}}$,  where $\rho$ is density, $v$ is velocity, and  $E$ is the specific energy of the system. The pressure~$p$  is given by the equation of state
$p=(\gamma-1)(E-\rho v^2/2)$, 
where $\gamma$~is the adiabatic constant that will be taken as $\gamma =1.4$. 
We now consider the interaction with a Mach~3 shock and a sine wave. The spatial domain is now given by $\Omega:=(-5,5)$, with the initial condition 
\begin{align*} 
(\rho,v,p) (x, 0) = 
\begin{cases}
  (27/7, 
    4\sqrt{35}/9, 
     31/3 ) & \text{if  $x\leq-4$,}  \\ 
   (1+ \sin(5x)/5, 0, 1) 
    & \text{if $x>-4$,} 
\end{cases} \end{align*} 
with left inflow and right outflow boundary conditions. This problem was first considered by Shu and Osher \cite{shuosher89}.

We run the simulation until $T=1.8$ and compare the schemes against a reference solution computed with a resolution of $N=16000$. Figures \ref{shuosher_o5_n=200}~(a) to~(d)  and~(e) correspond to resolutions of $N=200$ and $N=400$ points, respectively.
Both WENO-Z, YC-WENO and OWENO schemes produce  similar resolutions, being the one presented by the OWENO scheme slightly higher. The lowest resolution clearly corresponds to the JS-WENO scheme, especially for the case $N=200$. For $N=400$  the OWENO5 scheme appears to  capture  the shock  slightly better than the other schemes.

Finally, we show in Figure~\ref{shuosher_o5_n=200}~(c)   a comparison involving the error of each scheme with respect to the corresponding CPU time required to achieve it. We can see that the efficiency of all schemes 
is nearly the same in the case of fifth-order accuracy, although minor differences 
are found for lower resolution in benefit of both YC-WENO and OWENO schemes. 
Such asymptotic behaviour is probably due to the fact that there is no zero of order higher than one  along the derivative of the composition of the flux with the solution. All the  schemes considered can cope with the phenomena properly.

\subsubsection*{Example 8: Double Mach reflection problem}
We consider a test problem for the 2D Euler equations:
$$\boldsymbol{u}_t+\boldsymbol{f}^1(\boldsymbol{u})_x+\boldsymbol{f}^2(\boldsymbol{u})_y=0,$$
 with  \begin{align*} 
& \boldsymbol{u} =  \begin{pmatrix} 
    \rho \\
    \rho v^x \\
    \rho v^y \\
    E \end{pmatrix}, \quad 
  \boldsymbol{f}^1 (\boldsymbol{u})= \begin{pmatrix}  
    \rho v^x \\
    p+\rho (v^x)^2 \\
    \rho v^xv^y \\
    v^x(E+p) 
  \end{pmatrix},  \quad   \boldsymbol{f}^2 (\boldsymbol{u})= \begin{pmatrix} 
    \rho v^y \\
    \rho v^xv^y \\
    p+\rho (v^y)^2 \\
    v^y(E+p) 
  \end{pmatrix},
\end{align*}
where  $\rho$ is density, $(v^x, v^y)$  is velocity, $E$ is the specific energy, and  $p$  is pressure. The equation of state is
 $$p=(\gamma-1)\left(E-\frac{1}{2}\rho((v^x)^2+(v^y)^2)\right),$$ with  $\gamma =1.4$.

\begin{figure}[t]
  \centering
  \setlength{\tabcolsep}{1cm}
  \begin{tabular}{cc}
    \includegraphics[width=0.38\textwidth]{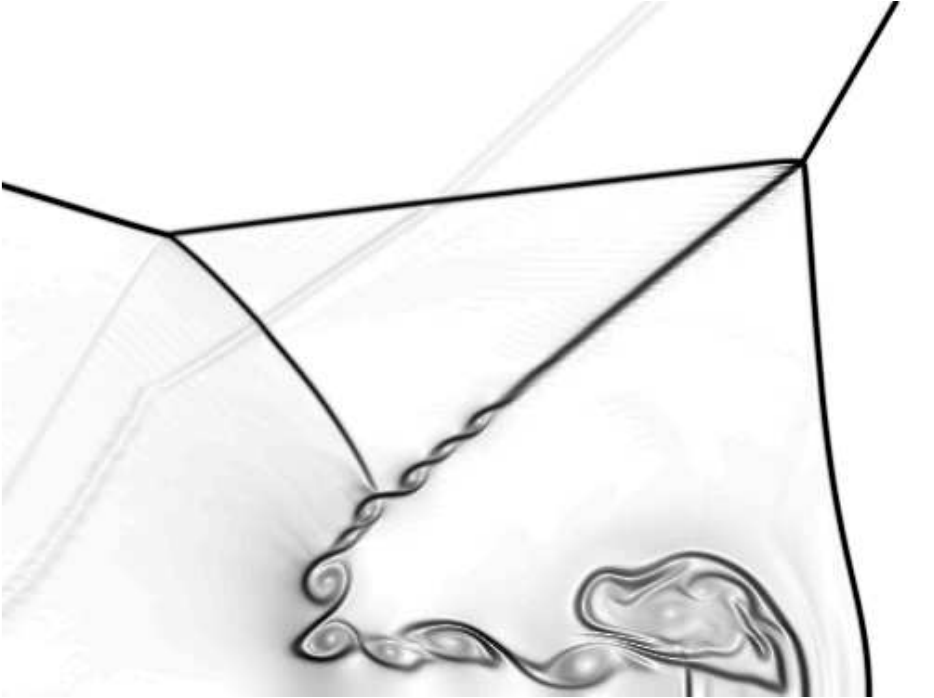} & \includegraphics[width=0.38\textwidth]{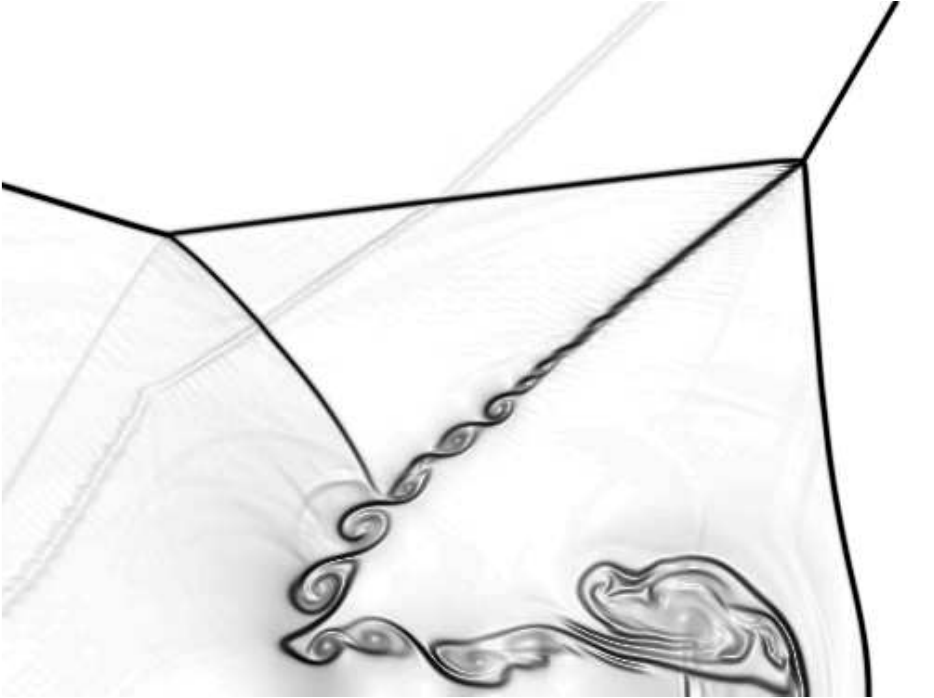} \\
    (a) JS-WENO5 & (b) WENO-Z5 \\
    \includegraphics[width=0.38\textwidth]{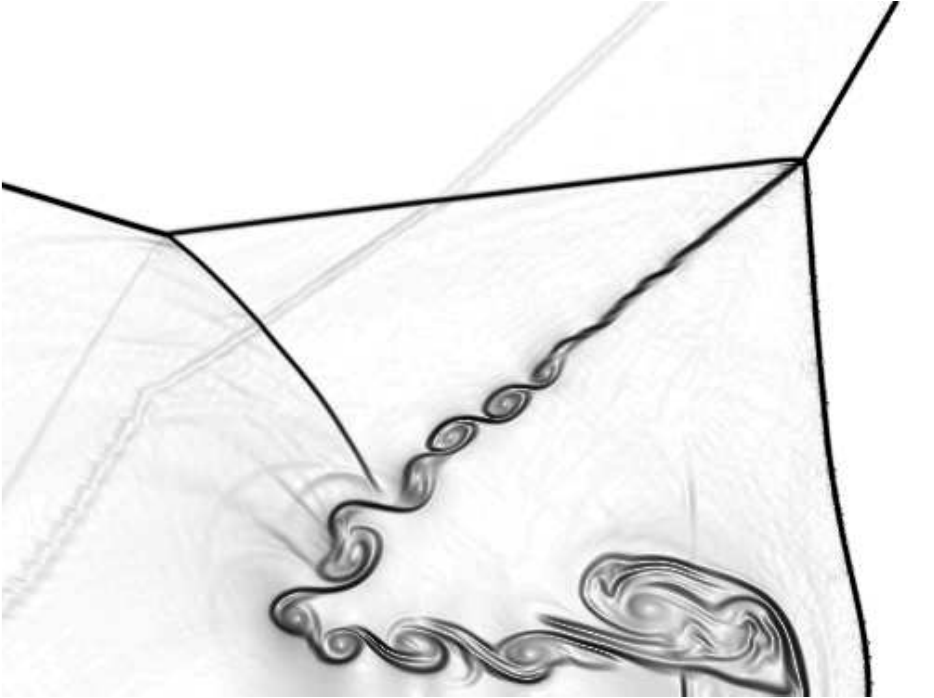} & \includegraphics[width=0.38\textwidth]{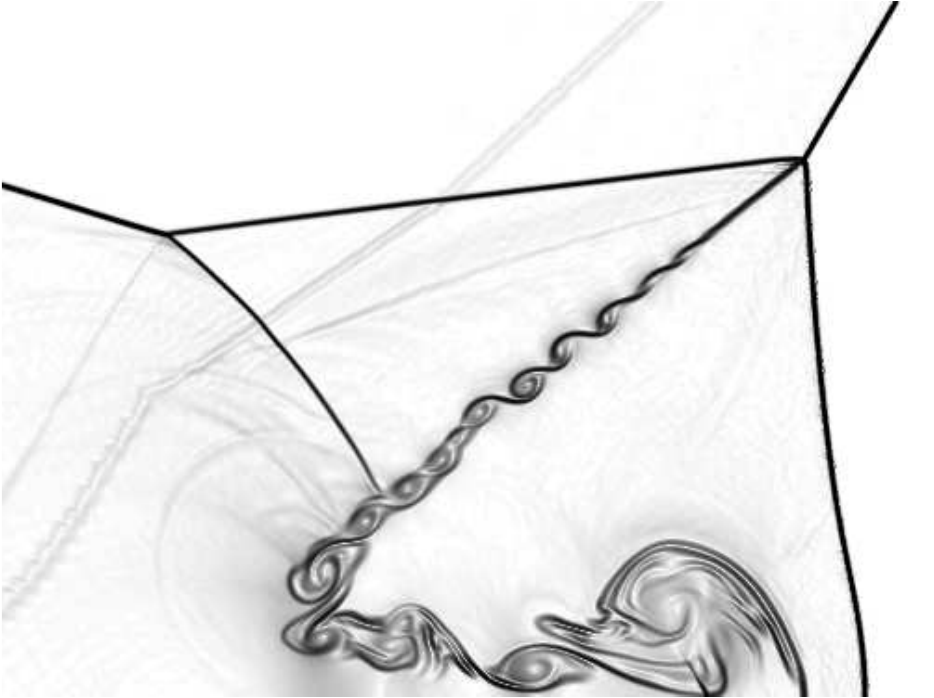} \\
    (c) YC-WENO5 & (d) OWENO5
  \end{tabular}
  \caption{Example~8 (Double Mach reflection problem, $2560\times640$, 2D Euler equations of gas dynamics): enlarged views  of the  turbulent zone of the numerical solutions at $T=0.2$ (Schlieren plot).}
  \label{dmr640}
\end{figure}

The Double Mach reflection test models a vertical right-going Mach 10 shock that hits an equilateral  triangle. By symmetry, we consider the problem defined only on the upper half part of the domain, which represents  a collision of the shock with a ramp with a slope of $30^{\circ}$ with respect to the horizontal line. Moreover,
we consider the equivalent problem defined in a rectangle but with the shock rotated $30^{\circ}$. The domain is the rectangle $\Omega=[0,4]\times[0,1]$, and the initial  conditions are given by
\begin{gather*}  (\rho,v^x,v^y,E) (x,y,0)=\begin{cases}
  \boldsymbol{c}_1= (\rho_1,v_1^x,v_1^y,E_1)    & \text{if $y\leq 1/4 +\tan(\pi/6)x$,} \\
  \boldsymbol{c}_2=  (\rho_2,v_2^x,v_2^y,E_2)     & \text{if $y > 1/4 +\tan(\pi/6)x$,}  
\end{cases} \\
    \boldsymbol{c}_1 =
    \bigl(8,8.25\cos(\pi/6),-8.25\sin(\pi/6),563.5\bigr), \quad 
    \boldsymbol{c}_2=    (1.4,0,0,2.5).
  \end{gather*} 
We impose inflow boundary conditions, with value $\boldsymbol{c}_1$, at the left side, $\{0\}\times[0,1]$, outflow boundary conditions both at $[0,1/4]\times\{0\}$ and $\{4\}\times[0,1]$, reflecting boundary conditions at  $(1/4,4]\times\{0\}$ and inflow boundary conditions at the upper side, $[0,4]\times\{1\}$, which mimics the shock at its actual traveling speed:
\begin{align*} 
 (\rho,v^x,v^y,E) (x,1,t)=\begin{cases}
  \boldsymbol{c}_1 & \text{if $x\leq 1/4 + (1+20t)/\sqrt{3}$,}  \\
  \boldsymbol{c}_2 & \text{if $x>1/4 + (1+20t)/\sqrt{3}$.}  
\end{cases} \end{align*}
We perform the simulations up to $T=0.2$  for the fifth order versions of 
JS-WENO, WENO-Z, YC-WENO method and our OWENO scheme, at a resolution of $2560\times640$ points, with results shown in 
Figure~\ref{dmr640}. A value $CFL=0.4$ has been used in all simulations.
The results show that WENO-Z, YC-WENO and OWENO schemes produce sharper resolution than JS-WENO, with OWENO presenting a slightly higher resolution with respect to YC-WENO, and in turn YC-WENO presenting a slightly higher resolution than WENO-Z. Table ~\ref{dmr_cpu} shows the CPU cost of the four schemes for the resolution of $128\times32$ points, in which it can be seen that the cost of all the involved schemes is similar.

\begin{table}[t] 
  \centering
  \begin{tabular}{|c|c|c|c|}
    \hline
    JS-WENO5 & WENO-Z5 & YC-WENO5 & OWENO5 \\
    \hline
    32.894029 & 34.199013 & 35.690326 & 36.847610 \\
    \hline
  \end{tabular}
  \smallskip 
  \caption{Example 8 (Double Mach reflection problem, $128\times32$, 2D Euler equations of gas dynamics): CPU cost comparison (in seconds).}
  \label{dmr_cpu}
\end{table}

\section{Conclusions}
\label{sec:conclusions}

We   propose novel WENO reconstructions, 
called OWENO reconstructions, in which the
accuracy is optimal regardless of the order of the critical
point to which the stencil converges. The approach is related to the
work by  Yamaleev and  Carpenter \cite{YamaleevCarpenter2009},  
We provide the necessary theoretical background to justify the properties of the scheme,
which outperforms related existing methods under some circumstances,
 both for smooth and discontinuous solutions, and behave similarly under other
  situations. 
  The fact that the new method does not always outperform existing ones 
   is consistent with the conclusions drawn 
  in \cite{BorgesCarmonaCostaEtAl2008}, where it is claimed that improvements in the numerical solution mainly
  depend on how far from zero are the weights associated to stencils crossed by
  discontinuities, rather than to the detection of critical points (especially if they are high-order critical points).
     However, this work finally presents a WENO reconstruction procedure which never 
  	  loses accuracy near critical points regardless of their order, relying only on the local 
  	  data and without any influence of scaling parameters such as tuning the 
  	  parameter~$\varepsilon$. Therefore, it closes the question of the maximal order
  	  that can be attained near critical points by means of WENO reconstructions.
  	   Some questions remain open, as for example the influence of the exponents $s_1$ and $s_2$ in the numerical
  	    dissipation and the determination their optimal values so as to reduce it as much as possible without generating artifacts or spurious oscillations.

       Nevertheless, we expect a much more significant improvement for third-order schemes, whose original version proposed by    Jiang and Shu \cite{JiangShu96}  loses order near first-order critical points, which in this case, unlike higher-order critical points, 
is a very common phenomenon appearing in solutions of any type of ordinary differential equations (ODEs) or PDEs. Therefore, fixing this issue would entail a substantial improvement in the case of third-order WENO schemes. Since the procedure that we have described here is not valid for the case of third-order schemes, we are currently working on the development of a third order scheme with unconditionally optimal accuracy for smooth data.

\appendix

\section{Technical results}\label{sec:appendix} The following results are 
 necessary for the development of the theoretical results
 presented in the main text, but their proofs  being quite technical and involved, have been postponed to this appendix to enhance the readability of the main text.

The following result, whose proof follows by 
  using Taylor expansion, is the key to proving Lemma~\ref{zr}.

\begin{lemma}\label{lemma2.2} 
 If
$
    \L\colon C^{m+1}[a, b] \to \Pi_{n}
 $
  is a linear and continuous operator with respect to 
  $\norm{\cdot}=\norm{\cdot}_{\infty}$, then there exists $K>0$ such that
  for any $\zeta\in[a,b]$ and $w \in [a,b]$, 
  \begin{align*}
    \L[f] (w) =\sum_{s=0}^{m}\frac{f^{(s)}(\zeta)}{s!}\L[(w-\zeta)^s] + \Delta_{m+1,\zeta}\L[f] 
     \quad \text{\em  with $\norm{\Delta_{m+1,\zeta}\L[f]}  \leq K \norm{\smash{f^{(m+1)}}}$.}  
    \end{align*} 
\end{lemma}

\begin{lemma}\label{zr}
  Let $a_0<a_1<\dots<a_n$ and  $z$ be fixed real numbers. Let
  $S=\{x_{0,h},\ldots,x_{n,h}\}$   be an $(n+1)$-point stencil 
  with $x_{i,h}=z+a_ih$ for $h>0$.  
  For any real function~$f$, assume that the reconstruction polynomial $p_h=p_h[f] \in \Pi_n$  satisfies  either
  $p_h(x_{i,h})=f(x_{i,h})$  for $i=0, \dots, n$ or 
  $$\int_{x_{i,h}-h/2}^{x_{i,h}+h/2}p_h(x)\,
  \mathrm{d}x=\int_{x_{i,h}-h/2}^{x_{i,h}+h/2}f(x)\,
  \mathrm{d}x \quad \text{\em for $i=0, \dots, n$,}$$
  depending on whether the data are point values \eqref{fpointvalues}
  or cell averages \eqref{fcellaverages}.
  Then, for $1\leq j \leq n$ and $s\geq j$, there exist polynomials $b_{s,j} \in
  \Pi_{n-j}$, depending uniquely on the type of reconstruction and
  parameters $a_0,\dots,a_n$, such that for any $f\in C^{m+1}$
  \begin{align} \label{lemma3.2ident} 
  p_h^{(j)}(z+wh)=\sum_{s=j}^mb_{s,j}(w)h^{s-j}f^{(s)}(z)+\bigO(h^{m+1-j})
\end{align}
for sufficiently small $wh$. The functions    $b_{s,j}$ have   the
following properties:   \begin{align*} b_{s,j}(w)=s!\binom{s}{j}w^{s-j} \quad \text{\em for $j\leq s\leq n$}, \end{align*}  
    and  $b_{s,1} \equiv 0$ if and only if $n=1$, $s$~is even and
      $a_0=-a_1$,     and $b_{s,1} \not \equiv 0$ otherwise.     
  \end{lemma}

\begin{proof}

 We  let $a=a_0-1/2$ and $b=a_n+1/2$ and define the 
operators
\begin{align*} 
\tilde{\L}_{\nu}, \L_{\nu,j}\colon C^{m+1}[a, b] \to \Pi_{n}, \quad  
\nu=1,2,  \quad j\geq 1
\end{align*}  
through the following  
conditions, where  $i=0,\dots,n$ and $j\leq n$: 
\begin{align} 		\label{theconds} 
  \tilde{\L}_1[f](a_i)&=f(a_i),\quad &\L_{1,j}[f]= \bigl(\tilde{\L}_{1}[f] \bigr)^{(j)},\\
  \int_{a_{i}-1/2}^{a_{i}+1/2}\tilde{\L}_2 [f](x) \, \mathrm{d}x &=
  \int_{a_{i}-1/2}^{a_{i}+1/2}f(x)\, \mathrm{d}x,\quad &\L_{2,j}[f]=\bigl(\tilde{\L}_{2}[f]
  \bigr)^{(j)}.    \label{theconds2} 
  \end{align}

The linearity of $\smash{\tilde{\L}}_{\nu}$ and $\L_{\nu,j}$ is clear and the continuity can be proven
by exploiting  conditions \eqref{theconds} and \eqref{theconds2}, 
e.g.,  by using 
Lagrange basis polynomials~$\varphi_i$ (standard ones for point evaluation); i.e.,  if we define 
  $\smash{\tilde\L_1[f] := f(a_0) \varphi_0+ \dots + f(a_n) \varphi_n}$,  
  then 
  \begin{align*} 
  \L_{1,j}[f]&=\sum_{i=0}^{n} f(a_i) \varphi_i^{(j)}, \quad 
  \bigl\Vert \L_{1,j}[f] \bigr\Vert \leq \max_{0\leq i\leq n}(f(a_i)) \sum_{i=0}^{n}
  \bigl\Vert\varphi_i^{(j)} \bigr\Vert\leq \Vert f\Vert \sum_{i=0}^{n}
  \bigl\Vert\varphi_i^{(j)} \bigr\Vert.
\end{align*}  
Similar arguments apply to the cell-average case ($\nu=2$).

With the notation $  S_{z,h}(w):=z+wh$, the polynomials  \eqref{lemma3.2ident} can be 
 expressed as $\smash{p_{h}=\tilde{\L}[f\circ S_{z,h}]\circ S_{z,h}^{-1}}$, 
  which means that 
    $\smash{p_{h}(x)=\tilde{\L}[f\circ S_{z,h}] ( (x-z)/h )}$,  
where either $\smash{\tilde\L}= \smash{\tilde\L}_1$ or
$\smash{\tilde\L} = \smash{\tilde\L}_2$, and correspondingly, either  $\L_j = \L_{1,j}$ or~$\L_j =\L_{2,j}$.  
Since  $\smash{(f\circ S_{z,h})^{(s)}(w)=h^{s}f^{(s)}(z+wh)}$, 
 Lemma~\ref{lemma2.2} for $\zeta=0$ yields 
\begin{align*}p_{h}^{(j)}(x)&=h^{-j}\tilde{\L}[f\circ S_{z,h}]^{(j)}  \bigl( (x-z)/h \bigr) =h^{-j}{\L_j}[f\circ S_{z,h}]  \bigl( (x-z)/h \bigr),\\
  p_{h}^{(j)}(z+wh)&=h^{-j}{\L_j}[f\circ S_{z,h}](w) \\ &  =
  h^{-j}\sum_{s=0}^{m}\frac{(f\circ S_{z,h})^{(s)}(0)}{s!}\L_j[w^s]
   +
  h^{-j}\Delta_{m+1,0}\L_j[f\circ S_{z,h}] \\ & =
  \sum_{s=j}^{m}h^{s-j}\frac{f^{(s)}(z)}{s!}\L_j[w^s] +\bigO(h^{m+1-j}),
\end{align*}
since $\smash{\tilde \L[w^s]=w^s}$ for $s\leq n$, therefore 
$\smash{\L_j[w^s]=(\tilde \L[w^s])^{(j)}=0}$ for $s<j$, and
$$\Delta_{m+1,0}\L_j[f\circ
S_{z,h}]\leq K \bigl\Vert (f\circ S_{z,h})^{(m+1)} \bigr\Vert_{[a,b]} =Kh^{m+1}
\bigl\Vert f^{(m+1)} \bigr\Vert_{S_{z,h}([a,b])}.$$
Therefore, the  result follows with $b_{s,j}(w)={\L_j[w^s]}/{s!}.$

Finally, if $n\geq 1$ and $b_{s,1}(w)=0$, then  for the first operator 
 we have 
 \begin{align*} 
 \tilde{\L}_1 [w^s]=\alpha \sii a_i^s=\alpha, \; i=0,\dots,n
 \sii  \text{$n=1$, $s$ is even and $a_0=-\alpha^{1/s}$, $a_1=\alpha^{1/s}$}.
 \end{align*} 
For the second operator, we have   $b_{s,1}(w)=0 \sii$
\begin{align*} 
 \tilde{\L}_2 [w^s]=\alpha   \sii
q_s(a_i) = (a_{i}+ 1/2 )^{s+1}- (a_{i}-  1/2 )^{s+1}=(s+1)\alpha, \quad 
  i=0,\dots,n, 
  \end{align*} 
where we define 
  \begin{align*}
q_s(x) & := (x+ 1/ 2)^{s+1}- (x- 1/2)^{s+1} = 
  \sum_{l=0}^{\lfloor s/2\rfloor}\binom{s+1}{2l+1}\frac{1}{2^{2l}}x^{s-2l}.
\end{align*}
Thus, 
by Rolle's
theorem,  there exist numbers $\smash{\tilde{a}}_i \in (a_{i-1}, a_{i})$,
$i=1,\dots,n$ such that $q_s'(\smash{\tilde{a}}_i)=0$. But 
\begin{align*}
  q_s'(x)&=
  \sum_{l=0}^{\lfloor s/2\rfloor}\binom{s+1}{2l+1}\frac{1}{2^{2l}}(s-2l)x^{s-2l-1}
\end{align*}
has only even-degree terms, with strictly positive coefficients, when
$s$ is odd (and therefore no roots) and only odd-degree terms, with strictly positive coefficients, when
$s$ is even (and therefore 0 as only root). This implies that  $s$ is
even, $n=1$ and 
$\smash{\tilde{a}}_1=0$, which yields $a_0 < \smash{\tilde{a}}_1=0 < a_1$. Since 
$q_s$ is an even function  and strictly increasing in $(0, \infty)$, for even $s$,
$q_s(a_0)=q_s(-a_0)=q_s(a_1)$ implies $a_1=-a_0$.  The converse
 is clear, since $n=1$, $a_1=-a_0$ and even $s$ implies that
$q_s(a_1)=q_s(a_0)=\alpha$ and therefore $\smash{\tilde{\L}}_2[w^s]=\alpha$
and $b_{s,1}(w)= (1 / s!)  \smash{\tilde{\L}}_2[w^s]'=0$.
\end{proof}

After some straightforward algebra,
we prove in the next result  that   $\bar\omega_i=\lim_{\eps\to 0}
  \omega_{i}$ exists and we obtain its rate of convergence.

  \begin{lemma}\label{lemma:pep1}
For  fixed data $f_{-r+1},\dots,f_{r-1}$,  we have $\omega_i = \bar\omega_i+\bigO(\eps^{s_2})$ and 
\begin{align}
  \label{eq:425}
  \bar\omega_i&=
  \begin{cases}
    c_i & \text{\em if $d_1d_2=0$,} \\
 \displaystyle  \frac{ c_{i}}{ \displaystyle \sum_{j=0, I_j\neq 0}^{r-1}c_j } & \text{\em if $d_1d_2\neq 0$,  $\exists k$
      with $I_k=0$, and $I_i= 0$,}     \\
    0& \text{\em if $d_1d_2\neq  0$,  $\exists k$
      with $I_k=0$, and $I_i\neq 0$,}     \\ \displaystyle
 \frac{c_i\big(  1 +
  \bar d/{I_{i}^{s_1}}\big)^{s_2}}{\displaystyle \sum_{j=0}^{r-1} c_j
\big(  1 +
\bar d/{I_{j}^{s_1}}\big)^{s_2}}    & \text{\em if $d_1d_2\neq 0$ and $I_k\neq 0$ for $k=0,\dots,r-1$.}
  \end{cases}
\end{align}

\end{lemma}

\begin{lemma}\label{lemma:1}
  If $f\in C^{s}(z)$ and 
   $\smash{f^{(s')}(z)=0}$  for all $s'<s$,  
   then
  \begin{align}
    e_i(h):=f(z+h/2)-p_i(z+h/2)&=\bigO(h^{\max \{ r,s \}}), \label{eihest}  \\ 
    e(h):=f(z+h/2)-q(z+h/2)&=\bigO(h^{\max \{r,s \}}).  \label{ehest} 
  \end{align}    
\end{lemma}
\begin{proof}
  We prove the result for the interpolatory case, the cell-average
  case is similar. Without loss of generality assume $z=0$.
  Using the Newton representation of the interpolation error, we get 
  \begin{align*}
    e_i(h)=f(x_{1/2})-p_i(x_{1/2}) 
     = \frac{f^{(r)}(\xi)}{r!} h^r    \prod_{l=0}^{r-1} \left(\frac{1}{2} - i + l 
     \right), 
  \end{align*}
  where $|\xi -z|
    < \max\{r-1-i, i\}h < rh$.
  The result follows for $s\leq r$. 
  For $s>r$, due to the assumption  and  using Taylor's remainder theorem, we get
  \begin{align*}
    f^{(r)}(\xi) = \frac{f^{(s)}(\xi_{s,r})}{(s-r)!}(\xi-z)^{s-r}
    |\xi_{s,r}-z|< |\xi-z|.
  \end{align*}
  It follows that for sufficiently small~$h_0$, 
  \begin{align*}
    \bigl|e_i(h) \bigr|\leq \max_{|\xi-z| < rh_0} \bigl|f^{(s)}(\xi) \bigr|\frac{r^s}{r!(s-r)!}h^s
     \quad \text{for $0<h<h_0$.} 
  \end{align*}
  This concludes the proof of \eqref{eihest}, 
   and \eqref{ehest} 
   follows from $\bar\omega_0 + \dots + \bar\omega_{r-1}=1$.
\end{proof}

In order to use the previous results, we consider
$x_{i,h}=z+(\alpha+i)h$, with $\alpha\in \mathbb R$ fixed and $i\in\mathbb Q$,
so that, for instance $x_{1/2,h}=z+(\alpha+1/2)h$. The reconstruction
polynomial $p_{r,i}$ associated to the  substencil
$S_{r,i}$ (see \eqref{Sridef})  corresponds to $p_h$ in Lemma
\ref{zr} for $n=r-1$ and 
\begin{align} \label{cijdef} 
a_j = a_{j,i} :=\alpha-r+i+1+j, \quad j=0,\dots,r-1. 
\end{align}

\begin{lemma}\label{roots_discont}
  Let $x_0<x_1<\dots<x_{n}$ be a stencil. Let $0\leq i_0\leq n-1$ and 
   $p \in \Pi_n$ be an interpolating polynomial  such that $p(x_i)=f_{\mathrm{L}}$ if $i\leq i_0$ and $p(x_i)=f_{\mathrm{R}}$ if $i>i_0$, with $f_{\mathrm{L}}\neq f_{\mathrm{R}}$. Then, $p^{(s)}$ has exactly $n-s$ roots, for $1\leq s\leq n$, and $p^{(s)} \in \bar{\Pi}_{n-s}$ for $0\leq s\leq n$. In particular, the parabola $p^{(n-2)}$ has two simple roots.
\end{lemma}

\begin{proof}
  Let $0\leq i\leq n-1$ such that $i\neq i_0$. Then, by construction,
  we have $p(x_i)=p(x_{i+1})$, and therefore by Rolle's theorem exists
  $\xi_i\in(x_i,x_{i+1})$ such that $p'(\xi_i)=0$, $0\leq i\leq
  n-1$. Therefore, $p' \in \Pi_{n-1}$ has at least $n-1$
  roots. However, since  $p$ takes different values it is not a
  constant polynomial, and   thus $p' \not\equiv 0$. Hence, $p' \in
  \bar{\Pi}_{n-1}$, $p'$ 
  must  have exactly $n-1$ roots and, a fortiori, $p \in
  \bar{\Pi}_{n}$.  A recursive application of
  Rolle's theorem yields 
  that $(p')^{(s-1)}=p^{(s)}\in\bar{\Pi}_{n-1-(s-1)}=\bar{\Pi}_{n-s}$ has exactly $(n-1)-(s-1)=n-s$ roots for
  $1\leq s\leq n$. 
\end{proof}

\begin{lemma}\label{linear}
  Let $x_{i,h}=z+a_ih$, $0\leq i\leq n$, be a grid with  $a_0<a_1<\dots<a_n$ and $p_h \in \Pi_n$ the interpolating polynomial such that $p_h(x_{i,h})=f_{i}$, for $f_{i}\in\R$, $0\leq i\leq n$. Then, given $0\leq s\leq n$, the $s$-th derivative of $P_h(w):=p_h(z+wh)$ can be written as
  $$P_h^{(s)}(w)=\sum_{j=0}^{n-s}L_{ \boldsymbol{a}}^{s,j}(f_{0,h},\ldots,f_{n,h})w^j, \quad 
   \boldsymbol{a} := (a_0,\ldots,a_n),$$
  with $L_{ \boldsymbol{a}}^{s,j}:\R^{n+1}\to\R$ a linear function, which does not depend on $h$.
  Furthermore,
  \begin{equation}\label{eq:Lrel}
    L_{ \boldsymbol{a} }^{s,j}(f_{0,h},\ldots,f_{n,h})=\frac{(s+j)!}{j!}L_{ \boldsymbol{a} }^{0,s+j}(f_{0,h},\ldots,f_{n,h}).
  \end{equation}
  Moreover, if $f_{i}=f(x_{i,h})$, for some $f\in C^{n+1}$, then
  \begin{align*} 
    L_{ \boldsymbol{a} }^{s,j}(f_{0,h},\ldots,f_{n,h})=\frac{h^{s+j}}{j!}f^{(s+j)}(z)+\bigO(h^{n+1}).
  \end{align*}
\end{lemma}

\begin{proof}
  Let $\mathcal{F}$ be the vector space of real functions and
    $\Phi_{\boldsymbol{a}}\colon \mathcal{F} \to \R^{n+1}$ 
  be the linear function given by
  $\Phi_{\boldsymbol{a}}(f)=(f(a_0),\dots,f(a_n))$. 
  Since $\ker  \Phi_{\boldsymbol{a}} \cap \Pi_n = 0$, and
  $\dim \Pi_n=n+1$, 
  $\Phi_{\boldsymbol{a}} |_{\Pi_n}$ is a bijection and
  $P_h=(\Phi_{\boldsymbol{a}} |_{\Pi_n})^{-1}(f_{0,h},\dots,f_{n,h})$. Since
  $\pi_i\colon \Pi_n \to \R$, $\pi_i(\sum_{j=0}^{n}\alpha_j
  w^j)=\alpha_i$ is a linear function, $\pi_i\circ
  (\Phi_{\boldsymbol{a}} |_{\Pi_n})^{-1}$ is also a linear function,
  therefore
  \begin{equation*}
    P_h(w)=\sum_{j=0}^{n}L_{\boldsymbol{a}}^{0,j}(f_{0,h},\dots,f_{n,h})w^j,\quad
    L_{\boldsymbol{a}}^{0,j} =\pi_i\circ
  (\Phi_{\boldsymbol{a}} |_{\Pi_n})^{-1}, 
\end{equation*}
from where 
equation \eqref{eq:Lrel} follows immediately.

Assume $f_{i}=f(x_{i,h})$, $f\in C^{n+1}(z)$. Since
$p_h(x)=P_h((x-z)/h)$,
  \begin{equation*}
    p_h(x)=\sum_{j=0}^{n}L_{\boldsymbol{a}}^{0,j}(f_{0,h},\dots,f_{n,h})h^{-j}(x-z)^j.
  \end{equation*}
This yields
  $L_{\boldsymbol{a}}^{0,j}(f_{0,h},\dots,f_{n,h})h^{-j}
  j!=p_h^{(j)}(z)$, for $j=0,\dots,n$.
   On the other hand the interpolation property yields
   $p_h^{(j)}(z)=f^{(j)}(z) + \bigO(h^{n+1-j})$, for $j=0,\dots,n$, thus implying
  \begin{equation*}
  L_{\boldsymbol{a}}^{0,j}(f_{0,h},\dots,f_{n,h})=\frac{f^{(j)}(z)}{j!}h^{j}
  + \bigO(h^{n+1}),
\end{equation*}
which, together with \eqref{eq:Lrel}, concludes the proof.
\end{proof}

\section*{Acknowledgements} 

AB, PM and DZ are supported by Spanish MINECO project
MTM2017-83942-P.  
RB is supported by  CRHIAM, project CONICYT/FON\-DAP/15130015;  CONICYT/PIA/AFB170001; and Fondecyt project 1170473. 
PM is also
supported by Conicyt (Chile),   pro\-ject PAI-MEC, folio 80150006.  
 DZ is also supported by Conicyt (Chile) through Fondecyt project 3170077.

\end{document}